\documentclass[preprint,12pt,authoryear]{elsarticle}

\journal{Journal of Number Theory}

\usepackage{amssymb}
\usepackage{amsfonts}

\usepackage[english]{babel}
\usepackage{amsmath,latexsym,enumerate,verbatim, amsthm}
\usepackage{stmaryrd}

\newcommand{\be}{\begin{enumerate}}
\newcommand{\ee}{\end{enumerate}}
\let\ds=\displaystyle

\def\N{{\mathbb N}} \def\Z{{\mathbb Z}}

\def\R{{\mathbb R}} \def\C{{\mathbb C}}

\def\lcm{\operatorname{lcm}}

\def\T{{\cal T}}

\def\s{{\bf s}}

\def\m{{\bf m}}
\def\i{{\bf i}}

\def\z{{\bf z}}
\def\e{{\bf e}}
\def\x{{\bf x}}
\def\y{{\bf y}}

\def\cb{{\bf c}}

\def\u{{\bf u}}

\def\c{{\bf c}}

\def\m{{\bf m}}
\def\X{{\bf X}}

\newcommand{\zerob}{\boldsymbol{0}}
\newcommand{\unb}{\boldsymbol{1}}

\newcommand{\alphab}{\boldsymbol{\alpha}}
\newcommand{\betab}{\boldsymbol{\beta}}

\newcommand{\nub}{\boldsymbol{\nu}}
\newcommand{\sigmab}{\boldsymbol{\sigma}}
\newcommand{\taub}{\boldsymbol{\tau}}
\newcommand{\lambdab}{\boldsymbol{\lambda}}

\def\eps{{ \varepsilon }}

\newcommand{\la}{\langle}

\newcommand{\ra}{\rangle}

\newtheorem{remark} {Remark}

\newtheorem{lemma}{Lemma}

\newtheorem{corollary}{Corollary}

\newtheorem{theorem}{Theorem}

\newtheorem{definition}{Definition}

\begin{document}
\begin{frontmatter}
%%***********************************************************
\title{Mean values of multivariable multiplicative functions and applications to the average number of cyclic subgroups and multivariable averages associated with the LCM function.}
\author[add1]{D. Essouabri }
\ead{driss.essouabri@univ-st-etienne.fr} 
\author[add1]{C. Salinas Zavala}
\ead{christoper.salinas.zavala@univ-st-etienne.fr}
\author[add2]{ L. T\'oth}
\ead{ ltoth@gamma.ttk.pte.hu}
%\date{}
%\maketitle
%\author{}
\address[add1]{Universit\'e Jean Monnet (Saint-Etienne), Institut Camille Jordan (UMR 5208 du CNRS), Facult\'e des Sciences et Techniques, 23 rue du Doc. Paul Michelon, Saint-Etienne, F-42023, France }
\address[add2]{Department of Mathematics, University of P\'ecs, Ifj\'us\'ag \'utja 6, P\'ecs, H-7624, Hungary }
%\affiliation[A]{organization={Univ. Lyon, UJM-Saint-Etienne, CNRS, Institut Camille Jordan UMR 5208},%Department and Organization
%addressline={Facult\'e des Sciences et Techniques, 23 rue du Doc. Paul Michelon}, 
        %    city={Saint-Etienne},
        %    postcode={F-42023}, 
       %     state={},
        %    country={France}}
            
% \affiliation[B]{organization={Department of Mathematics, University of P\'ecs},%Department and Organization
        %    addressline={Ifj\'us\'ag \'utja 6}, 
        %    city={P\'ecs},
       %     postcode={H-7624}, 
         %   state={},
       %     country={Hungary}}

\begin{abstract}
We use multiple zeta functions to prove, under suitable assumptions,  precise asymptotic formulas for the averages of multivariable multiplicative functions. As applications, we prove some conjectures on the average number of cyclic subgroups of the group $\Z_{m_1}\times\dots\times \Z_{m_n}$ and multivariable averages associated with  the LCM function.
\end{abstract}
\end{frontmatter}

%\vskip .1 in
\setcounter{tocdepth}{2}

\section{Introduction}
Our paper is motivated by the following recent results and conjectures. Let $n\in \N$ and for $m_1,\dots, m_n \in \N$ let $c_n(m_1, \dots , m_n) $  denote 
the number of cyclic subgroups of the group $\Z_{m_1}\times \dots \times \Z_{m_n}$. W. G. Nowak and L. T\'oth \cite{NT} proved the asymptotic formula
$$\sum_{1\leq m_1,m_2\leq x} c_2(m_1,m_2) =x^2\left(\frac{12}{\pi^4} (\ln x)^3+a_2  (\ln x)^2+a_1  (\ln x)+a_0\right) +O(x^{\frac{1117}{701}+\eps}) 
\quad {\mbox{ as }}\quad x\rightarrow \infty,$$
where $a_0,a_1$ and $a_2$ are explicit constants.  This error term was improved by L. T\'oth and W. Zhai \cite{TZ2018} into $O(x^{3/2} (\ln x)^{13/2})$. 
The case $n=3$ was investigated by L. T\'oth and W. Zhai \cite{TZ2020} showing that 
$$\sum_{1\leq m_1,m_2,m_3\leq x} c_3(m_1,m_2,m_3) =x^3 \sum_{j=0}^7 c_j (\ln x)^j+ O(x^{8/3+\varepsilon}),$$
where $c_j$ ($0\le j \le 7$) are explicit constants. For the proof they used a multidimensional Perron formula and the complex integration method. 
It is natural to conjecture that such a result holds for $n\geq 4$. 

T. Hilberdink, F. Luca, and L. T\'oth \cite{HLT} investigated the following three averages associated with the LCM function:
\begin{equation}\label{snx}
S_n(x):=\sum_{1\leq m_1,\dots, m_n \leq x} \frac{1}{\lcm(m_1,\dots,m_n)},
\end{equation}
\begin{equation}\label{unx}
U_n(x):=\sum_{1\leq m_1,\dots, m_n \leq x\atop \gcd(m_1,\dots,m_n)=1} \frac{1}{\lcm(m_1,\dots,m_n)},
\end{equation}
and
\begin{equation}\label{vnx}
V_n(x):=\sum_{1\leq m_1,\dots, m_n \leq x} \frac{m_1\dots m_n}{\lcm(m_1,\dots,m_n)}.
\end{equation}
By using the convolution method, they obtained in their paper asymptotic formulas with error terms for $S_2(x)$, $U_2(x)$ and $V_2(x)$.
For $n\geq 3$, they only obtained the estimates
$$(\ln x)^{2^n-1}\ll S_n(x)\ll (\ln x)^{2^n-1}, \quad (\ln x)^{2^n-2}\ll U_n(x)\ll (\ln x)^{2^n-2}, 
$$ 
$$
x^n\ll V_n(x)\ll x^n~(\ln x)^{2^n-2} \quad {\mbox{ as }}\quad x\rightarrow \infty,$$
and conjectured that asymptotic formulas with error terms also exist for these three averages for $n\geq 3$.

In order to prove these conjectures, we introduce a reasonably large class of mutivariable multiplicative functions (see Definition \ref{class-def}). 
For a function $f: \N^n\to \R_+$ in this class, we establish in Theorem \ref{main1} the existence of the meromorphic continuation of the associated  multiple zeta function  
$$\s=(s_1,\dots,s_n)\to \mathcal M (f;\s):=\sum_{m_1\geq 1,\dots,m_n\geq 1}\frac{f(m_1,\dots,m_n)}{m_1^{s_1}\dots m_n^{s_n}}$$
and derive several precise properties of this meromorphic continuation. By combining our Theorem \ref{main1} and La Bret\`eche's  multivariable Tauberian Theorem (i.e., Theorems 1 and 2 of \cite{bretechecompo}) {\it we deduce in our Theorem 2 a precise asymptotic formula for the multivariable average 
$$\mathcal N_\infty (f; x):=\sum_{\m=(m_1,\dots,m_n)\in \N^n\atop \|\m\|_\infty=\max_i m_i \leq x} f(m_1,\dots,m_n) \quad {\mbox{as }} x\rightarrow \infty,$$ and 
derive from it four corollaries.}

{\it Our first application, namely Corollary \ref{corollary1}, establishes the conjecture concerning the number of cyclic subgroups of the group $\Z_{m_1}\times \dots \times \Z_{m_n}$, in any dimension $n$. Our Corollaries \ref{corollary2}, \ref{corollary3} and \ref{corollary4} prove the conjectures on the three sums above associated with the LCM
 function}.

Variants of Theorem \ref{main2} with other norm choices can be obtained by combining our Theorem \ref{main1} and the first author's  multivariable tauberian theorem (i.e., Corollary 2 of \cite{Toricmanin}). For example, for the class of H\"older's norms
$\|\x\|_d:=\sqrt[d]{|x_1|^d+\dots+|x_n|^d}$ ($d\geq 1$), {\it we obtain in Theorem \ref{main3} an asymptotic for the multivariable average 
$$\mathcal N_d (f; x):=\sum_{\m=(m_1,\dots,m_n)\in \N^n\atop \|\m\|_d =\sqrt[d]{m_1^d+\dots+m_n^d}\leq x}
f(m_1,\dots,m_n)\quad {\mbox{as }} x\rightarrow \infty.$$
As an application of Theorem \ref{main3}, we derive in Corollaries  \ref{corollary5} and  \ref{corollary6} the analogues of Corollaries \ref{corollary1} and \ref{corollary4} for the H\"older's norms $\|~\|_d $}.

\subsection{Notations}
\be
\item
$\N=\{1,2,\dots \}$,
$\N_0=\N \cup \{0\}$; $\R_+ = [0, \infty).$
\item
The expression: $ f(\lambda,{\bf y},{\bf x}){\ll}_{{}_{{\bf y}}} g({\bf x})$
uniformly in ${ \bf x}\in X $ and ${\lambda}\in \Lambda$
means there exists $A=A({\bf y})>0$,
such that,
$\forall {\bf x}\in X {\mbox { and }}\forall {\lambda}\in
{\Lambda}\quad |f(\lambda,{\bf y},{\bf x})|\leq Ag({\bf x}) $;
\item Let $d\in[1,+\infty[$, for any ${\bf x}=(x_1,...,x_n) \in \R^n$, we set
$\|{\bf x}\|_d=\sqrt[d]{|x_1|^d+...+|x_n|^d}$,
 and $\|{\bf x}\|_\infty =\max_{i=1,\dots,n} |x_i|$.
  We denote the canonical basis of $\R^n$ by
$(\e_1,\dots,\e_n)$ (i.e. $\e_{i,j} =1$ if $i=j$ and $\e_{i,j} =0$ if $i\neq j$). The standard inner product on $\R^n$ is
denoted by $\la .,. \ra$. We set also
$\zerob=(0,\dots,0)$ and $\unb =(1,\dots,1)$;
\item
We denote a vector in $\C^n$ by $\s=(s_1,\dots,s_n)$, and write
 $\s={\sigmab}+i{\taub},$
 where ${\sigmab } = (\sigma_1,\dots,\sigma_n)$ and
 ${\taub }=(\tau_1,\dots,\tau_n)$ are the real resp. imaginary
 components of $\s$ (i.e. $\sigma_i=\Re(s_i)$ and
$\tau_i=\Im(s_i)$ for all $i$). We also write
$\la \x, \s\ra$ for $\sum_i x_i s_i$  if $\x \in \R^n, \s \in \C^n$;
\item A function $f: \N^{n} \rightarrow \C$ is said to be
  multiplicative if for all $\m=(m_1,\dots,m_n)\in \N^n$ and
  $\m'=(m_1',\dots,m_n')\in \N^n$ satisfying
$\gcd\left(\lcm\left(m_i\right),\lcm\left(m_i'\right)\right)=1$ we have
$f\left(m_1 m_1',\dots , m_n m_n'\right)=f(\m) \cdot f(\m')$;

\item Let $F$ be a meromorphic function on a domain ${\mathcal D}$ of
  $\C^n$ and let ${\mathcal S}$ be the support of its polar divisor.
$F$ is said to be of moderate growth if
$\exists a,b>0$ such that $\forall \delta >0$,
 $F(\s) \ll_{\sigmab,\delta} 1+\|\taub\|_1^{a\|\sigmab\|_1+b}$
uniformly in $\s=\sigmab+i\taub\in {\mathcal D}$
verifying $d(\s, {\mathcal S})\geq \delta $;
\ee

\section{A class of multivariable multiplicative functions and statement of the main results}
\subsection{A class of multivariable multiplicative functions}
To simplify the exposition, we introduce first the following three definitions.
\begin{definition}\label{data-def}
A quadruple $(g,\kappa, \c, \delta)$ is said to be a data if
\be
\item $g : \N_0^n\to \N_0$ is a function of subexponential growth; that is $g$ verifies for any $\eps >0$ $g(\nub)\ll_\eps
 e^{\eps \|\nub\|_1}$ uniformly in $\nub\in \N_0^n$;
\item  $\kappa : \N_0^n \to [1,\infty)\cup\{0\}$ is a function verifying $\kappa(\zerob)=0$ and
$\inf_{\nub \in \N_0^n \setminus\{\zerob\}} \frac{\kappa(\nub)}{\|\nub\|_1} >0$
\item $\c=(c_1,\dots,c_n)\in [0,\infty)^n$ and  $\delta \in (0,\infty)$.
\ee
\end{definition}
We now introduce the class of multivariable multiplicative functions on which we will focus in this paper.\\
\begin{definition}\label{class-def}
Let $(g,\kappa, \c, \delta)$ be a data as in definition \ref{data-def}.\\
A multivariable multiplicative function $f:\N^n\to \R$ is said to be in the class $\mathcal C (g,\kappa, \c, \delta)$ if
for any $\eps >0$,
\begin{equation}\label{mainassumption}
f(p^{\nu_1},\dots,p^{\nu_n}) - g(\nub)~p^{\la \c, \nub\ra -\kappa (\nub)}\ll_{\eps} e^{\eps \|\nub\|_1} ~p^{\la \c, \nub\ra -\kappa (\nub) -\delta},
\end{equation}
uniformly in $\nub \in \N_0^n$ and $p$ prime number.
\end{definition}
We will need also the following  integral definition.
\begin{definition}\label{Infxdef}
Let $I$ be a finite subset of $\N_0^n\setminus \{\zerob\}$, $\u=\left(u(\nub)\right)_{\nub \in I}$ be a finite sequence of positive integers and $\c=(c_1,\dots,c_n)\in [0,\infty)^n$.
We denote by $\nub^1,\dots,\nub^r$ the elements of $I$ where $r=\# I$,
and  define  the finite sequence $q_k$ $(0\leq k\leq r)$ by
$$
q_0=0 \quad {\mbox{ and }} \quad q_k=\sum_{j=1}^k u(\nub^j) ~\forall k=1,\dots, r.
$$
We define then for $x>0$ the integral
$$
\mathcal I_n(I,\u, \c; x):=\int_{\mathcal A(I,\u;x)} \frac{dy_1\dots dy_{q_r}}{ \prod_{k=1}^r \prod_{\ell =q_{k-1}+1}^{q_k}
y_\ell^{1-\la \nub^k, \c\ra}},
$$
where $\ds \mathcal A (I,\u;x):=\left\{\y\in [1,\infty)^{q_r};~~ \prod_{k=1}^r \prod_{\ell =q_{k-1}+1}^{q_k} y_\ell^{\la \nub^k, \e_j\ra} \leq x ~~\forall j=1,\dots,n\right\}$.
\end{definition}
\subsection{Statement of the main results}
Let $f:\N^n\to \R$ be a multivariable multiplicative function. We assume that $f$ belongs to the class  $\mathcal C (g,\kappa, \c, \delta)$ associated to the data $(g,\kappa, \c, \delta)$ (see definitions 1 and 2 above).\\
We assume also that the finite set
\begin{equation}\label{Inonempty}
I=I(\kappa, g):=\{\nub \in \N_0^n \mid \kappa(\nub)=1{\mbox{ and }} g(\nub)\neq 0\} \quad {\mbox {is nonempty.}}
\end{equation}

The following theorem is {\it the main analytic ingredient} of this paper:
\begin{theorem}\label{main1}
\be
\item the multiple zeta function
$$\s=(s_1,\dots,s_n)\to \mathcal M (f;\s):=\sum_{m_1\geq 1,\dots,m_n\geq 1}\frac{f(m_1,\dots,m_n)}{m_1^{s_1}\dots m_n^{s_n}}$$
converges absolutely in the domain $\{\s\in \C^n \mid \Re (s_i) > c_i ~\forall i=1,\dots,n\}$;
\item there exists $\eps_0>0$ such that the function $$\s=(s_1,\dots,s_n)\to \mathcal H(f,\c;\s):=\left(\prod_{\nub \in I} \la \nub, \s\ra^{g(\nub)}\right)~\mathcal M (f;\c+\s)$$
has {\bf holomorphic} continuation to the domain $\{\s\in \C^n \mid \Re (s_i) > -\eps_0 ~\forall i=1,\dots,n\}$ and verifies in it the following estimate: for all $\eps >0$,
$$\mathcal H(f,\c;\s) \ll_\eps \prod_{\nub \in I} \left(|\la \nub, \s\ra|+1\right)^{g(\nub)\left(1-\frac{1}{2}\min\left(0, \Re (\la \nub, \s\ra)\right)\right)+\eps};$$
\item $\mathcal H(f,\c;\zerob)$ is given by the following convergent Euler product:
\begin{equation}\label{c0}
\mathcal H(f,\c;\zerob)=\prod_p \left(1-\frac{1}{p}\right)^{\sum_{\nub \in I} g(\nub)} \left(\sum_{\nub \in \N_0^n}
\frac{f(p^{\nu_1},\dots,p^{\nu_n})}{p^{\la \nub, \c\ra}}\right).
\end{equation}
\ee
\end{theorem}

Combining our Theorem \ref{main1} and La Bret\`eche's  multivariable Tauberian Theorem (i.e Theorems 1 and 2 of \cite{bretechecompo}) yields to the following multivariable mean value theorem:
\begin{theorem}\label{main2}
Let $f:\N^n\to \R_+$ be a nonnegative multivariable multiplicative function satisfying assumptions of Theorem \ref{main1}.
Set $J:=\{\e_i \mid  c_i=0\}$ where $(\e_1,\dots,\e_n)$ is the canonical basis of $\R^n$.
Set also $\rho:=\left(\sum_{\nub \in I} g(\nub)\right) +\#J -Rank\left(I\cup J\right)$.\\
Then,
there exist a polynomial $Q_\infty$ of degree at most $\rho$ and a positive constant $\mu_\infty >0$ such that
$$\mathcal N_\infty (f; x):=\sum_{\m=(m_1,\dots,m_n)\in \N^n\atop \|\m\|_\infty=\max_i m_i \leq x} f(m_1,\dots,m_n)
=x^{\|\c\|_1} Q_\infty(\ln x) + O\left(x^{\|\c\|_1-\mu_\infty}\right)\quad {\mbox{ as }} x\rightarrow \infty.$$

Furthermore, if we assume in addition that the two following assumptions hold:
\be
\item $Rank\left(I\cup J\right) =n$;
\item $\unb=(1,\dots,1)$ is in the interior of the cone generated by $I\cup J$; that is $ \unb \in con^*\left(I\cup J\right):=\{\sum_{\nub \in I\cup J}\lambda_{\nub} \nub \mid \lambda_{\nub} \in (0,\infty)~\forall \nub \in I\cup J\}$,
\ee
Then, the degree of the polynomial $Q_\infty$ is equal to $\rho=\left(\sum_{\nub \in I} g(\nub)\right) +\#J -n$ and the main term of $\mathcal N_\infty (f; x)$ is given by
$$\mathcal N_\infty (f; x) = C_n(f)  K_n(f,\| \|_\infty) ~x^{\|\c\|_1} (\ln x)^{\rho}+ O\left((\ln x)^{\rho-1}\right) \quad {\mbox{as }} x\rightarrow \infty,$$
where $C_n(f):= \mathcal H(f,\c;\zerob)>0$ is defined by the Euler product (\ref{c0}) and
$$K_n(f,\| \|_\infty):=\lim_{x\rightarrow \infty} {\mathcal I}_n(I,\u,\c; x) ~x^{-\|\c\|_1} (\ln x)^{-\rho} >0,\qquad {\mbox{ where}}$$
$\mathcal I_n(I,\u,\c; x)$ is the integral (see definition \ref{Infxdef}) associated to the finite set $I$, the finite sequence
$\u=\left(g(\nub)\right)_{\nub \in I}$ and to the vector $\c$.
\end{theorem}

\begin{remark}\label{assumptionssymp}
The existence of the limit $K_n(f,\| \|_\infty)$ follows from the proof of Theorem\ref{main2}.\\
If $\{\e_1,\dots,\e_n\} \subset I\cup J$, then the two assumptions $Rank\left(I\cup J\right) =n$ and $\unb \in con^*\left(I\cup J\right)$ clearly hold.
\end{remark}\par

Combining our Theorem \ref{main1} and the first author's  multivariable tauberian theorem (i.e corollary 2 of \cite{Toricmanin}) yields to the following multivariable mean value theorem for H\"older's norms
$\|\x\|_d:=\sqrt[d]{|x_1|^d+\dots+|x_n|^d}$ ($d\geq 1$):
\begin{theorem}\label{main3}
Let $f:\N^n\to \R_+$ be a nonnegative multivariable multiplicative function satisfying assumptions of Theorem \ref{main1}.\\
Assume that $\c=(c_1,\dots,c_n)\in (0,\infty)^n$.
Set
\be
\item
$\rho:=\left(\sum_{\nub \in I} g(\nub)\right) - Rank (I)$;
\item $I_{\cb}:=\{\la\c, \nub\ra^{-1} \nub \mid \nub \in I\}$ and $\u:=\left(u(\betab)\right)_{\betab \in I_{\cb}}$ where
$\ds u(\betab)=\sum_{\nub \in I, \la\c, \nub\ra^{-1} \nub =\betab} g(\nub)$.
\ee
Then,
there exist a polynomial $Q$ of degree at most $\rho$ and a positive constant $\mu >0$ such that
$$\mathcal N_d (f; x):=\sum_{\m=(m_1,\dots,m_n)\in \N^n\atop \|\m\|_d =\sqrt[d]{m_1^d+\dots+m_n^d}\leq x}
f(m_1,\dots,m_n)
=x^{\|\c\|_1} Q(\ln x) + O\left(x^{\|\c\|_1-\mu}\right)\quad {\mbox{ as }} x\rightarrow \infty.$$

Furthermore, if we assume in addition that $Rank ( I )=n$ and $ \unb \in con^*\left(I\right)$,
then, the degree of the polynomial $Q$ is equal to $\rho=\left(\sum_{\nub \in I} g(\nub)\right) -n$ and the main term of $\mathcal N_d (f; x)$ is given by
$$\mathcal N_d (f; x) = C_n(f)  K_n(f,\| \|_d) ~x^{\|\c\|_1} (\ln x)^{\rho}+ O\left((\ln x)^{\rho-1}\right) \quad {\mbox{as }} x\rightarrow \infty,$$
where $C_n(f):= \mathcal H(f,\c;\zerob)>0$ is defined by the Euler product (\ref{c0}) above and
$$  K_n(f,\| \|_d):=\left(\prod_{\nub \in I} \la \nub, \cb\ra^{-g(\nub)}\right)~\frac{d^{\rho+1}~
 A_0(\T_{\cb} ,P_d)}{\|\c\|_1 ~ \rho!} > 0. $$
 where  $A_0(\T_{\cb}, P_d)>0$ is the mixed volume constant (see \S  2.3.3  of \cite{Toricmanin}) associated to the pair $\mathcal T_{\cb} :=(I_{\cb}, \u)$ and the polynomial 
 $P_d= X_1^{d}+\dots + X_n^{d}$.
\end{theorem}
\subsection{Applications}
We will now give the applications that motivated our general results of section \S 2.2.
\subsubsection{On the  average number of cyclic subgroups of the group $\Z_{m_1}\times\dots\times \Z_{m_n}$}
Let $n\in \N$. For $m_1,\dots, m_n \in \N$ denote by $c_n(m_1, \dots , m_n) $  the number of cyclic subgroups of the group $\Z_{m_1}\times \dots \times \Z_{m_n}$. Set $$ G_n(x):=\sum_{1\leq m_1,\dots, m_n \leq x}c_n(m_1,\dots,m_n).$$
As we mentioned in the introduction,  precise asymptotic for $G_2(x)$ was obtained by W. G. Nowak and L. T\'oth  in \cite{NT} and  improved by L. T\'oth and W. Zhai in \cite{TZ2018}. The case $n=3$ was also investigated by L. T\'oth and W. Zhai in \cite{TZ2020}. It is natural to conjecture that such a result holds for $n\geq 4$. The following result establish this conjecture in any dimension $n$.
\begin{corollary}\label{corollary1}
Let $n \in \N$. There exists a polynomial $Q_1$ of degree $2^n-1$ and $\mu_1 >0$ such that
$$G_n(x):=\sum_{1\leq m_1,\dots, m_n \leq x}c_n(m_1,\dots,m_n) = x^n ~Q_1(\ln x) + O(x^{n-\mu_1}) \quad {\mbox{ as }} \quad x\rightarrow \infty.$$
In particular, we have
$$ G_n(x) =C_n(c_n)  K_n(c_n,\| \|_\infty)~x^n (\ln x)^{2^n-1} +O\left( x^n (\ln x)^{2^n-2}\right) \quad {\mbox{ as }} \quad x\rightarrow \infty,$$
where
\begin{equation}\label{besoincn}
C_n(c_n):= \prod_p \left(1-\frac{1}{p}\right)^{2^n+n-1} \left(\sum_{\nub \in \N_0^n}
\frac{c_n(p^{\nu_1},\dots,p^{\nu_n})}{p^{\|\nub\|_1}}\right) >0
\end{equation}
 and
 $$K_n(c_n,\| \|_\infty):=\lim_{x\rightarrow \infty} {\mathcal I}_n(I,\u; x) ~x^{-n} (\ln x)^{-2^n+1} >0,\qquad {\mbox{ where}}
 $$
$\ds \mathcal I_n(I,\u, \c; x)$ is the integral (see definition \ref{Infxdef}) associated to $\ds I=\{0,1\}^n\setminus \{\zerob\}$,  to the sequence
$\ds \u=\left(u(\nub)\right)_{\nub \in I}$ defined by
$u(\e_i)=2$ $\forall i=1,\dots,n$ and $u(\nub) =1$ $\forall \nub \in I\setminus \{\e_1,\dots, \e_n\}$
and to the vector $\c=\unb$.
\end{corollary}
\begin{remark}\label{casen=2or3incor1}
We will compute more explicitly in \S 7 below the constants $C_n(c_n)$ and $K_n(c_n,\| \|_\infty)$ for $n=2$ and $n=3$. In particular, we will prove in \S 7.1 and \S 7.3 that $C_2(c_2)=\frac{36}{\pi^4}$ and $K_2(c_2,\| \|_\infty)=\frac{1}{3}$. Thus, our mains term in the asymptotic of $G_2(x)$ agree with the main term obtained by the convolution method in \cite{NT} by W. G. Nowak and L. T\'oth.
\end{remark}

\subsubsection{Some multivariable averages associated to the LCM function}
As we mentioned in the introduction,  T. Hilberdink, F. Luca, and L. T\'oth introduced in \cite{HLT} the three averages (\ref{snx}), (\ref{unx}) and (\ref{vnx}) associated to the LCM function and obtained in this paper asymptotic formulas for $S_2(x)$, $U_2(x)$ and $V_2(x)$.
For $n\geq 3$, they only obtained the following estimates
$$(\ln x)^{2^n-1}\ll S_n(x)\ll (\ln x)^{2^n-1}, \quad (\ln x)^{2^n-2}\ll U_n(x)\ll (\ln x)^{2^n-2}, 
$$
$$
x^n\ll V_n(x)\ll x^n~(\ln x)^{2^n-2},$$
and conjectured that asymptotic formulas also exist for these three averages for $n\geq 3$.\\
The following three corollaries prove these conjectures.
\begin{corollary}\label{corollary2}
Let $n\in \N$. There exists a polynomial $Q_2$ of degree $2^n-1$ and $\mu_2 >0$ such that
$$S_n(x):=\sum_{1\leq m_1,\dots, m_n \leq x} \frac{1}{\lcm(m_1,\dots,m_n)} = Q_2(\ln x) + O(x^{-\mu_2}) \quad {\mbox{ as }} \quad x\rightarrow \infty.$$
In particular, we have
$$ S_n(x) =C_n(s_n)  K_n(s_n,\| \|_\infty)~(\ln x)^{2^n-1} +O\left( (\ln x)^{2^n-2}\right) \quad {\mbox{ as }} \quad x\rightarrow \infty,$$
where
\begin{equation}\label{besoinsn}
C_n(s_n):=\prod_p \left(1-\frac{1}{p}\right)^{2^n-1} \left(\sum_{k=0}^\infty
\frac{(k+1)^n -k^n}{p^{k}}\right)>0,
\end{equation}
 and
 $$K_n(s_n,\| \|_\infty):=\lim_{x\rightarrow \infty} {\mathcal I}_n(I,\u,\c; x) ~(\ln x)^{-2^n+1} >0,\qquad {\mbox{ where}}$$
$\ds \mathcal I_n(I,\u, \c; x)$ is the integral (see definition \ref{Infxdef}) associated to $\ds I=\{0,1\}^n\setminus \{\zerob\}$,  to the sequence
$\ds \u=\left(u(\nub)\right)_{\nub \in I}$ defined by
$u(\nub) =1$ $\forall \nub \in I$
and to the vector $\c=\zerob$.
\end{corollary}

\begin{corollary}\label{corollary3}
Let $n \in \N\setminus\{1\}$. There exists a polynomial $Q_3$ of degree $2^n-2$ and $\mu_3 >0$ such that
$$U_n(x):=\sum_{1\leq m_1,\dots, m_n \leq x\atop \gcd(m_1,\dots,m_n)=1} \frac{1}{\lcm(m_1,\dots,m_n)} = Q_3(\ln x) +
 O(x^{-\mu_3}) \quad {\mbox{ as }} \quad x\rightarrow \infty.$$
In particular, we have
$$ U_n(x) =C_n(u_n)  K_n(u_n,\| \|_\infty)~(\ln x)^{2^n-2} +O\left( (\ln x)^{2^n-3}\right) \quad {\mbox{ as }} \quad x\rightarrow \infty,$$
where
\begin{equation}\label{besoinun}
C_n(u_n):=\prod_p \left(1-\frac{1}{p}\right)^{2^n-1} \left(\sum_{k=0}^\infty
\frac{(k+1)^n -k^n}{p^{k}}\right)>0,
\end{equation}
 and
 $$K_n(u_n,\| \|_\infty):=\lim_{x\rightarrow \infty} {\mathcal I}_n(I,\u,\c; x) ~(\ln x)^{-2^n+2} >0,\qquad {\mbox{ where}}
 $$
$\ds \mathcal I_n(I,\u, \c; x)$ is the integral (see definition \ref{Infxdef}) associated to $\ds I=\{0,1\}^n\setminus \{\zerob, \unb\}$,  to the sequence
$\ds \u=\left(u(\nub)\right)_{\nub \in I}$ defined by
$u(\nub) =1$ $\forall \nub \in I$
and to the vector $\c=\zerob$.
\end{corollary}

\begin{corollary}\label{corollary4}
Let $n \in \N$. There exists a polynomial $Q_4$ of degree $2^n-n-1$ and $\mu_4 >0$ such that
$$V_n(x):=\sum_{1\leq m_1,\dots, m_n \leq x} \frac{m_1\dots m_n}{\lcm(m_1,\dots,m_n)} = x^n~Q_4(\ln x) + O(x^{n-\mu_4}) \quad {\mbox{ as }} \quad x\rightarrow \infty.$$
In particular, we have
$$ V_n(x) = C_n(v_n)  K_n(v_n,\| \|_\infty)~x^n~(\ln x)^{2^n-n-1} +O\left( x^n~(\ln x)^{2^n-n-2}\right) \quad {\mbox{ as }} \quad x\rightarrow \infty,$$
where
\begin{equation}\label{besoinvn}
C_n(v_n):=\prod_p \left(1-\frac{1}{p}\right)^{2^n-1} \left(\sum_{k=0}^\infty
\frac{(k+1)^n -k^n}{p^{k}}\right)>0,
\end{equation}
 and
 $$K_n(v_n,\| \|_\infty):=\lim_{x\rightarrow \infty} {\mathcal I}_n(I,\u,\c; x) ~x^{-n} (\ln x)^{-2^n+n+1} >0,
 \qquad {\mbox{ where}}
 $$
$\ds \mathcal I_n(I,\u, \c; x)$ is the integral (see definition \ref{Infxdef}) associated to $\ds I=\{0,1\}^n\setminus \{\zerob\}$,  to the sequence
$\ds \u=\left(u(\nub)\right)_{\nub \in I}$ defined by
$u(\nub) =1$ $\forall \nub \in I$
and to the vector $\c=\unb$.
\end{corollary}
\begin{remark}\label{casen=2or3incor2-4}
The constants  $C_n(s_n)$,  $C_n(u_n)$ and  $C_n(v_n)$ are equal.
We will compute more explicitly in sections  \S 7.2,  \S 7.4,  \S 7.5 and \S 7.6 below the constants $C_n(.)$ and $K_n(.,\| \|_\infty)$ for $n=2$ and $n=3$. More precisely, we will prove that
\be
\item $ C_2(s_2)=C_2(u_2)=C_2(v_2)=\frac{6}{\pi^2}$ and \\
$ C_3(s_3)=C_3(u_3)=C_3(v_3)=\prod_p\left(1-\frac{9}{p^2}+\frac{16}{p^3}-\frac{9}{p^4}+\frac{1}{p^6}\right)$;
\item $K_2(s_2,\| \|_\infty)=\frac{1}{3}$, $K_2(u_2,\| \|_\infty)=1$ and $K_2(v_2,\| \|_\infty)=1$;
\item $K_3(s_3,\| \|_\infty)=\frac{11}{3366}$, $K_3(u_3,\| \|_\infty)=\frac{11}{480}$ and $K_3(v_3,\| \|_\infty)=\frac{1}{16}$.
\ee
In particular, our mains terms in the asymptotic of $S_2(x)$, $U_2(x)$ and $V_2(x)$ agree with those obtained by the convolution method  in \cite{HLT} by T. Hilberdink, F. Luca, and L. T\'oth.
\end{remark}

\subsubsection{Multivariable averages with other norms}

The following two results give analogues of corollaries \ref{corollary1} and \ref{corollary4} for some other choices of norms.

\begin{corollary}\label{corollary5}
Let $n\in \N$ and $d\geq 1$.
There exists a polynomial $Q_5$ of degree $2^n-1$ and $\mu_5 >0$ such that
$$G_{n,d}(x):=\sum_{\m=(m_1,\dots,m_n)\in \N^n\atop \|\m\|_d =\sqrt[d]{m_1^d+\dots+m_n^d}\leq x} c_n(m_1,\dots,m_n) = x^n ~Q_5(\ln x) + O(x^{n-\mu_5}) \quad {\mbox{ as }} \quad x\rightarrow \infty.$$
Moreover, if we set
${\tilde{I}}:=\left\{ \|\nub\|_1^{-1}\nub;~\nub \in \{0,1\}^n\setminus\{\zerob\}\right\}$ and
$\u=\left(u(\betab)\right)_{\betab \in {\tilde{I}}}$ where $u(\betab)=2$ if $\betab \in \{\e_1,\dots,\e_n\}$ and  $u(\betab)=1$ otherwise,
then
$$ G_{n,d}(x) = C_n(c_n)~K_n(c_n,\| \|_d)~x^n (\ln x)^{2^n-1} +O\left( x^n (\ln x)^{2^n-2}\right) \quad {\mbox{ as }} \quad x\rightarrow \infty,$$
where $C_n(c_n)>0$ is given by (\ref{besoincn}) and
$$K_n(c_n,\| \|_d)=\left(\prod_{k=2}^n k^{-{n\choose k}} \right) \frac{~d^{2^n}~ A_0(\T,P_d)}{n~ (2^n-1)!}>0,$$
 where $A_0(\T,P_d)$ is the mixed volume constant (see \S  2.3.3  of \cite{Toricmanin}) associated to the pair $\mathcal T =(\tilde{I}, \u)$ and the polynomial
 $P_d= X_1^{d}+\dots + X_n^{d}$.
 \end{corollary}

\begin{corollary}\label{corollary6}
Let $n\in \N$ and $d\geq 1$.
There exists a polynomial $Q_6$ of degree $2^n-n-1$ and $\mu_6 >0$ such that
$$V_{n,d}(x):=\sum_{\m=(m_1,\dots,m_n)\in \N^n\atop \|\m\|_d =\sqrt[d]{m_1^d+\dots+m_n^d}\leq x}
\frac{m_1\dots m_n}{\lcm(m_1,\dots,m_n)} = x^n ~Q_6(\ln x) + O(x^{n-\mu_6}) \quad {\mbox{ as }} \quad x\rightarrow \infty.$$
Moreover, if we set
${\tilde{I}}:=\left\{ \|\nub\|_1^{-1}\nub;~\nub \in \{0,1\}^n\setminus\{\zerob\}\right\}$ and
$\u=\left(u(\betab)\right)_{\betab \in {\tilde{I}}}$ where $u(\betab)=1$ $\forall \betab \in {\tilde{I}}$,
then
$$ V_{n,d}(x) = C_n(v_n)~K_n(v_n,\| \|_d)~x^n (\ln x)^{2^n-n-1} +O\left( x^n (\ln x)^{2^n-n-2}\right) \quad {\mbox{ as }} \quad x\rightarrow \infty,$$
where $C_n(v_n)>0$ is given by (\ref{besoinvn}) and
$$K_n(v_n,\| \|_d)=\left(\prod_{k=2}^n k^{-{n\choose k}} \right)  \frac{~d^{2^n-n}~ A_0(\T,P_d)}{n~ (2^n-n-1)!} >0,$$
 where $A_0(\T,P_d)$ is the mixed volume constant (see \S  2.3.3  of \cite{Toricmanin}) associated to the pair $\mathcal T =(\tilde{I}, \u)$ and the polynomial
 $P_d= X_1^{d}+\dots + X_n^{d}$.
 \end{corollary}

 \begin{remark}\label{casen=2or3incor5-6}
The constants  $C_n(.)$,  are independent on the choice of the norm. The constants $K_n(.)$ depend on the choice of the norm. We will compute more explicitly in \S 7.7 below the constants $K_n(c_n,\|  \|_d)$ and $K_n(v_n,\|  \|_d)$ for $n=2$ and $n=3$. More precisely, we will prove that
\be
\item $ K_2(c_2,\| \|_d)=\dfrac{1}{6d}\dfrac{\Gamma\left(1/d\right)^2}{\Gamma\left(2/d\right)}$ and 
$K_2(v_2,\| \|_d) =\dfrac{1}{2d}\dfrac{\Gamma(1/d)^2}{\Gamma(2/d)}$;
\item $ K_3(c_3,\| \|_d)=\dfrac{31~\Gamma\left(1/d\right)^3}{30240~ d^{2}~\Gamma\left(3/d\right)}$ and
$K_3(v_3,\| \|_d)=\dfrac{\Gamma\left(1/d\right)^3}{2 ~d^2 ~\Gamma\left(3/d\right)}$.
\ee
\end{remark}

\section{Proof of Theorem \ref{main1}}
Let $f:\N^n\to \R$ be a multiplicative function in the class $\mathcal C (g,\kappa, \c, \delta)$. \\
Define for $\nub\in \N_0^n$ and $p$ prime $V(p,\nub)$ by the formula
\begin{equation}\label{vpnub-def}
f(p^{\nu_1},\dots,p^{\nu_n}) =\left( g(\nub)+ V(p,\nub)\right)~p^{\la \c, \nub\ra -\kappa (\nub)}
\end{equation}
Since  $(g,\kappa, \c, \delta)$ is a data, point 1 of definition \ref{data-def} and assumption (\ref{mainassumption}) can then be written in the following more convenient equivalent form:
\begin{equation}\label{mainassumptionbis}
\forall \eps>0, \quad g(\nub) \ll_{\eps} e^{\eps \|\nub\|_1} \quad {\mbox{ and }} \quad
V(p,\nub) \ll_{\eps} e^{\eps \|\nub\|_1} ~p^{-\delta},
\end{equation}
uniformly in $\nub \in \N_0^n$ and in $p$ prime number.

Moreover, point 2 of definition \ref{data-def} implies that there exists $\beta >0$ such that
\begin{equation}\label{kappanublowerbound}
\kappa (\nub) \geq \max \left( 1, \beta \|\nub\|_1 \right) \quad \forall \nub \in \N_0^n\setminus\{\zerob\}.
\end{equation}

\subsection{Proof of point 1 of Theorem \ref{main1}}
Let $\s=(s_1,\dots,s_n)\in \C^n$ be such that $\sigma_i=\Re(s_i)>c_i$ $\forall i=1,\dots,n$.\\
Set $\sigmab=(\sigma_1,\dots, \sigma_n)$ and $\eta=\frac{1}{2}\min_{i=1,\dots,n} (\sigma_i -c_i)>0$.\\
So, we have $\sigma_i \geq  c_i+2\eta$ $\forall i$ and $\la \sigmab, \nub\ra \ge \la \c, \nub\ra +2\eta \|\nub\|_1$ $\forall \nub \in \N_0^n$.
Choose $\eps>0$ small enough such that $e^{\eps}<2^{\eta}$.
It follows then from (\ref{vpnub-def}) and (\ref{mainassumptionbis}) that  we have for any  prime number $p$,
\begin{eqnarray*}
\sum_{p}\sum_{\|\nub\|_1 \geq 1}\left|\frac{f(p^{\nu_1},\dots,p^{\nu_n})}{p^{\la \s ,\nub\ra}}\right|
&= &\sum_{p}\sum_{\|\nub\|_1 \geq 1}\frac{\left|f(p^{\nu_1},\dots,p^{\nu_n})\right|}{p^{\la
    \sigmab, \nub\ra}}
\ll_\eps \sum_{p}\sum_{\|\nub\|_1 \geq 1}\frac{e^{\eps \|\nub\|_1} p^{\la \c, \nub\ra -\kappa (\nub)}}{p^{\la
    \cb, \nub\ra+2\eta \|\nub\|_1}}\\
    &\ll_\eps& \sum_{p}\sum_{\|\nub\|_1 \geq 1}\frac{e^{\eps \|\nub\|_1} }{p^{\kappa(\nub)+2\eta \|\nub\|_1}}
\ll_\eps  \sum_{p}\frac{1}{p^{1+\eta}}\sum_{\|\nub\|_1 \geq 1}\frac{e^{\eps \|\nub\|_1} }{p^{\eta \|\nub\|_1}}\\
&\ll_\eps & \sum_{p} \frac{1}{p^{1+\eta}}\sum_{\|\nub\|_1 \geq 1}\left(\frac{e^{\eps}}{2^{\eta}}\right)^{\|\nub\|_1}
\ll_\eps \sum_{p} \frac{1}{p^{1+\eta}}<\infty.
\end{eqnarray*}
The multiplicativity of $f$ implies then that $\s\to \mathcal M(f;\s)$ converges absolutely and that
\begin{equation}\label{eulerprod}
\mathcal M(f;\s)=\sum_{\m \in \N^{n}}
\frac{f(m_1,\dots,m_n)}{m^{s_1}\dots m_n^{s_n}}=
\prod_p \left(\sum_{\nub \in \N_0^n} \frac{f(p^{\nu_1},\dots,p^{\nu_n}) }{p^{\la \s, \nub\ra}}\right).
 \end{equation}
    This ends the proof of point 1 of Theorem \ref{main1}.\qed

\subsection{Two useful lemmas}
Recall that $I=I(\kappa, g):=\{\nub \in \N_0^n \mid \kappa(\nub)=1 {\mbox{ and }} g(\nub)\neq 0\}$ is an nonempty set.\\
For all $t \in \R$, set $U_t :=\{\s \in \C^n \mid \sigma_i=\Re(s_i) >t ~\forall i=1,\dots,n\}$.\\
We need the following two lemmas:
\begin{lemma}\label{lemma1}
There exists $\eps_1, \eta_1>0$ such that for any prime number $p$, the function
$$\s\mapsto R_p(\s):=\left(\sum_{\|\nub\|_1\geq 1} \frac{f(p^{\nu_1},\dots,p^{\nu_n}) }{p^{\la \c+\s, \nub\ra}}\right) -\left(\sum_{\nub \in I} \frac{g(\nub)}{p^{1+\la\nub, \s\ra}}\right)$$
is holomorphic in the domain $U_{-\eps_1}$ and verifies in it the estimate
$$ R_p(\s) \ll p^{-1-\eta_1} \quad {\mbox{ uniformly in }} p.$$
\end{lemma}

\begin{lemma}\label{lemma2}
Set $\eps_2=\inf_{\nub \in I} \frac{1}{4\|\nub\|_1}$ and $\eta_2=\frac{1}{2}$. Then, for any prime number $p$, the function
$$\s\mapsto L_p(\s):=\left(\prod_{\nub \in I} \left(1-\frac{1}{p^{1+\la\nub, \s\ra}}\right)^{g(\nub)}\right) -1 +
\left(\sum_{\nub \in I} \frac{g(\nub)}{p^{1+\la\nub, \s\ra}}\right)$$
is holomorphic in the domain $U_{-\eps_2}$ and verifies in it the estimate
$$ L_p(\s) \ll p^{-1-\eta_2} \quad {\mbox{ uniformly in }} p.$$
\end{lemma}
\subsubsection{Proof of Lemma \ref{lemma1}}
Fix $\beta >0$ such that  (\ref{kappanublowerbound}) holds. Fix also a positive  integer $N$ verifying
$\ds N\geq \max \left(4 \beta^{-1}, \max_{\nub \in I}\|\nub\|_1\right).$\\
Identity (\ref{vpnub-def}) implies that for $p$ prime number and $\s\in U_0=\{\s\in \C^n \mid \sigma_i >0 ~\forall i\}$, we have
\begin{equation}\label{decompo-Rp}
 R_p(\s) = R_p^1(\s)+R_p^2(\s), \quad {\mbox{ where}}
 \end{equation}
$$ R_p^1(\s)=\sum_{1\leq \|\nub\|_1\leq N } \frac{V(p,\nub)}{p^{\kappa(\nub)+\la\nub, \s\ra}}+\sum_{\nub \not \in I\atop 1\leq \|\nub\|_1\leq N} \frac{g(\nub)}{p^{\kappa(\nub)+\la\nub, \s\ra}}\quad
{\mbox{ and }}\quad
R_p^2(\s)=\sum_{\|\nub\|_1> N} \frac{g(\nub)+V(p,\nub)}{p^{\kappa(\nub)+\la\nub, \s\ra}}.$$
To prove Lemma \ref{lemma1} it suffices to verify that both $\s\mapsto R_p^1(\s)$ and $\s\mapsto R_p^2(\s)$ satisfy its  conclusions.

{\bf CLAIM 1:} $\s\mapsto R_p^1(\s)$ satisfies the conclusions of  Lemma \ref{lemma1}.\\
{\bf Proof of CLAIM 1:} It's clear that $\s \mapsto R_p^1(\s)$ is holomorphic in the whole space $\C^n$.\\
Let $\eps >0$. It follows from (\ref{mainassumptionbis}) and (\ref{kappanublowerbound}) that for $p$ prime number and $\s \in U_{-\eps}=\{\s\in \C^n\mid \sigma_i >-\eps ~\forall i\}$, we have
\begin{equation}\label{tata}
|R_p^1(\s)|\leq \sum_{1\leq \|\nub\|_1 \leq N} \frac{|V(p,\nub)|}{p^{1+\la\nub, \sigmab\ra}}+\sum_{\nub \not \in I\atop 1\leq \|\nub\|_1\leq N} \frac{g(\nub)}{p^{\kappa(\nub)+\la\nub, \sigmab\ra}}
\ll \sum_{\nub \in I} \frac{p^{-\delta}}{p^{1-\eps\|\nub\|_1}}+\sum_{\nub \not \in I,~g(\nub)\neq 0\atop 1\leq \|\nub\|_1\leq N}
\frac{1}{p^{\kappa(\nub)-\eps\|\nub\|_1}}
\end{equation}
Since $\kappa(\nub)>1$ if $\nub \not \in I\cup\{\zerob\}$ and $g(\nub)\neq 0$, it is clear that we can choose $\eps >0$ small enough such that
$$\mu_1:=\min_{\nub \in I} (\delta-\eps\|\nub\|_1) >0  {\mbox{ and }}
\mu_2:=\min\{\kappa(\nub)-\eps\|\nub\|_1 -1 \mid 1\leq \|\nub\|_1\leq N,~\nub \not \in I {\mbox { and }} g(\nub)\neq 0\}>0.$$
Set $\mu =\min (\mu_1,\mu_2)>0$. It follows then from (\ref{tata}) that we have
$R_p^1(\s) \ll p^{-1-\mu}$ uniformly in $p$ prime number and in $\s\in U_{-\eps}$. This ends the proof of CLAIM 1.\qed

{\bf CLAIM 2:} $\s\mapsto R_p^2(\s)$ satisfies the conclusions of  Lemma \ref{lemma1}.\\
{\bf Proof of CLAIM 2:}
Fix $\eps >0$ such that $e^\eps < 2^{\beta/4}$. Assumptions  (\ref{mainassumptionbis}) and (\ref{kappanublowerbound}) imply that we have uniformly in $p$ prime number and in
$\s\in U_{-\beta/2}$,
\begin{eqnarray*}
\sum_{\|\nub\|_1> N} \left|\frac{g(\nub)+V(p,\nub)}{p^{\kappa(\nub)+\la\nub, \s\ra}}\right|
&\ll_\eps& \sum_{\|\nub\|_1> N} \frac{e^{\eps \|\nub\|_1}}{p^{\beta \|\nub\|_1+\la\nub, \sigmab\ra}}
\leq \sum_{\|\nub\|_1> N} \frac{e^{\eps \|\nub\|_1}}{p^{\frac{\beta}{2} \|\nub\|_1}}
\ll_\eps \frac{1}{p^{\frac{\beta}{4} N}}\sum_{\|\nub\|_1> N} \frac{e^{\eps \|\nub\|_1}}{2^{\frac{\beta}{4} \|\nub\|_1}}\\
&\ll_\eps& \frac{1}{p^{\frac{\beta}{4} N}}\sum_{\|\nub\|_1> N} \left(\frac{e^{\eps}}{2^{\frac{\beta}{4}}}\right)^{\|\nub\|_1}
\ll_\eps \frac{1}{p^{\frac{\beta}{4} N}} \leq  \frac{1}{p^2}.
\end{eqnarray*}
We deduce that $\s\mapsto R_p^2(\s)$ is holomorphic in $ U_{-\beta/2}$ and verifies the estimates
$ R_p^2(\s) \ll p^{-2}$ uniformly in $p$ prime number and $\s \in U_{-\beta/2}$. This ends the proof of CLAIM 2 and also ends the proof of Lemma \ref{lemma1}.\qed

\subsection{Proof of Lemma \ref{lemma2}}
It is clear that $\s\mapsto L_p(\s)$ is holomorphic in $\C^n$ for any $p$. \\
Set now $\eps_2=\inf_{\nub \in I} \frac{1}{4\|\nub\|_1}$. It follows that  for $\s \in U_{-\eps_2} {\mbox{ and }} \nub \in I$,
$1+\la \nub , \sigmab \ra \geq 1-\eps_2 \|\nub\|_1 \geq 3/4.$\\
Newton Binomial theorem implies then that we have uniformly in
$\s \in U_{-\eps_2}$ and in $p$ prime number,
\begin{eqnarray*}
| L_p(\s)| &=&\left|\left(\prod_{\nub \in I} \left(1-\frac{1}{p^{1+\la\nub, \s\ra}}\right)^{g(\nub)}\right) -1 +
\left(\sum_{\nub \in I} \frac{g(\nub)}{p^{1+\la\nub, \s\ra}}\right)\right|\\
&=& \left|\sum_{0\leq k_{\nub} \leq g(\nub) ~\forall \nub \in I,\atop \sum_{\nub \in I} k_{\nub} \geq 2}
\frac{\prod_{\nub\in I} (-1)^{k_{\nub}}{g(\nub) \choose k_{\nub}}}
{p^{\sum_{\nub \in I} k_{\nub}(1+\la \nub , \s \ra)}}\right|
\ll \sum_{0\leq k_{\nub} \leq g(\nub) ~\forall \nub \in I,\atop \sum_{\nub \in I} k_{\nub} \geq 2}
\frac{1}{p^{\sum_{\nub \in I} k_{\nub}(1+\la \nub , \sigmab \ra)}}\\
&\ll& \sum_{0\leq k_{\nub} \leq g(\nub) ~\forall \nub \in I,\atop \sum_{\nub \in I} k_{\nub} \geq 2}
\frac{1}{p^{\frac{3}{4}\sum_{\nub \in I} k_{\nub}}} \ll \frac{1}{p^{3/2}}.
\end{eqnarray*}
This ends the proof of Lemma \ref{lemma2}.
\subsubsection{Proof of parts 2 and 3 of Theorem \ref{main1}}
Define the function $\s=(s_1,\dots,s_n)\mapsto \mathcal E(f;\s)$ by
\begin{equation}\label{efs-def}
\mathcal E(f;\s):=\left(\prod_{\nub \in I} \zeta \left(1+\la \nub, \s\ra\right)^{-g(\nub)}\right) ~\mathcal M(f;\c+\s).
\end{equation}
Part 1 of Theorem \ref{main1} implies then that $\s\mapsto \mathcal E(f;\s)$ converges absolutely in the domain $U_0=\{\s\in \C^n \mid \sigma_i >0 ~\forall i\}$. Moreover,  The multiplicativity of $f$  imply  that for all
$\s \in U_0 $:
\begin{equation}\label{eulerproddec1}
\mathcal E(f;\s) = \prod_p \mathcal E_p(f;\s), \quad {\mbox{where  }}
\end{equation}
$$\mathcal E_p(f;\s):=\prod_{\nub \in I} \left(1-\frac{1}{p^{1+\la \nub , \s \ra}}\right)^{g(\nub)}~
\left(\sum_{\nub \in \N_0^n} \frac{f(p^{\nu_1},\dots,p^{\nu_n})}{p^{\la \nub ,\cb+\s\ra}}\right).$$
We will now prove the following needed lemma:
\begin{lemma}\label{lemma3}
There exists $\eps_0>0$ such that the Euler product $\s\mapsto \mathcal E(f;\s) = \prod_p \mathcal E_p(f;\s)$ converges absolutely and defines a bounded holomorphic function in  the domain \\
$U_{-\eps_0}=\{\s\in \C^n\mid \sigma_i >-\eps_0~\forall i=1,\dots,n\}$.
\end{lemma}
{\bf Proof of Lemma \ref{lemma3}:}\\
We will use in the sequel of this proof notation of Lemmas 1 and 2.
Lemmas 1 and 2 imply that for any prime $p$ and any $\s\in U_0$,
\begin{eqnarray}\label{epfsdec2}
\mathcal E_p(f;\s)&=& \left(1-\left(\sum_{\nub \in I} \frac{g(\nub)}{p^{1+\la\nub, \s\ra}}\right)+L_p(\s)\right)
\left( 1+\left(\sum_{\nub \in I} \frac{g(\nub)}{p^{1+\la\nub, \s\ra}}\right)+ R_p(\s)\right)\nonumber\\
&=& 1-A_p(\s)^2+ B_p(\s)+C_p(\s), \quad  {\mbox{ where }}
\end{eqnarray}
$$
A_p(\s):=\sum_{\nub \in I} \frac{g(\nub)}{p^{1+\la\nub, \s\ra}},~
B_p(\s):= \left(1-A_p(\s)\right) R_p(\s)  {\mbox{ and }}
C_p(\s):= L_p(\s) \left( 1+A_p(\s)+ R_p(\s)\right).
$$
Let $\eps_1,\eps_2,\eta_1,\eta_2 >0$ the positive constants defined in Lemmas 1 and 2. \\
Set $\eps_0 =\min (\eps_1,\eps_2)>0$ and $\eta_0 =\min (\eta_1,\eta_2)=\min (\eta_1, 1/2)>0$.
Lemmas 1 and 2 imply that the three function $A_p$, $B_p$ and $C_p$ are holomorphic in $U_{-\eps_0}$ and that we have uniformly in $p$ prime number and $\s\in U_{-\eps_0}$  the following estimates:
\be
\item
$\ds A_p(\s)\ll \sum_{\nub \in I} \frac{1}{p^{1+\la\nub, \sigmab\ra}}\leq
\sum_{\nub \in I} \frac{1}{p^{1-\eps_0\|\nub\|_1}}\ll \frac{1}{p^{3/4}}$;
\item $\ds A_p(\s)^2\ll \frac{1}{p^{3/2}} \ll \frac{1}{p^{1+\eta_0}}$;
\item $\ds B_p(\s)\ll \left(1+ \frac{1}{p^{3/4}}\right) \frac{1}{p^{1+\eta_1}}\ll  \frac{1}{p^{1+\eta_0}}$;
\item $\ds C_p(\s)\ll \frac{1}{p^{1+\eta_2}} \left( 1+\frac{1}{p^{3/4}}+ \frac{1}{p^{1+\eta_1}}\right)\ll
\frac{1}{p^{1+\eta_0}}$.
\ee

It follows that for any prime number $p$, the function $\s\mapsto \mathcal E_p(f;\s) -1$ is holomorphic in $U_{-\eps_0}$ and
verifies
$\ds  \mathcal E_p(f;\s) -1 \ll \frac{1}{p^{1+\eta_0}} \quad {\mbox{ uniformly in }} \s \in U_{-\eps_0} {\mbox{ and in the prime number }} p.
$\\
We deduce that the Euler product $\s\mapsto \mathcal E(f;\s) = \prod_p \mathcal E_p(f;\s)$ converges absolutely and defines a bounded holomorphic function in  $U_{-\eps_0}$. This ends the proof of Lemma \ref{lemma3}.\qed \par
We are now ready to prove points 2 and 3 of Theorem \ref{main1}.
Combining part 1 of Theorem 1, (\ref{efs-def}) and (\ref{eulerproddec1}) implies that for $\s \in U_0$,
\begin{equation}\label{bienf}
\mathcal H(f,\c;\s):=\left(\prod_{\nub \in I} \la \nub, \s\ra^{g(\nub)}\right)~\mathcal M (f;\c+\s)\\
=\left(\prod_{\nub \in I} \left(\la \nub, \s\ra \zeta \left(1+\la \nub, \s\ra\right)\right)^{g(\nub)}\right)
~\mathcal E(f;\s).
\end{equation}
Part 2 of Theorem \ref{main1} follows then from Lemma \ref{lemma3} and the following two classical properties of Riemann zeta function:
$s\mapsto s\zeta(1+s)$ is holomorphic in $\C$ and verifies in the half-plane $\{\Re (s) >-1\}$ the estimate $s~\zeta(1+s)\ll_\eps (1+|s|)^{1-\frac{1}{2}\min\left(0,\Re (s)\right)+\eps}$, $\forall \eps>0$.\\
Moreover, since $s \zeta(1+s)|_{s=0}=1$, we deduce from (\ref{bienf}) and (\ref{eulerproddec1}) that
$$\mathcal H(f,\c;\zerob)=\mathcal E(f;\zerob)=
\prod_p  \left(1-\frac{1}{p}\right)^{\sum_{\nub \in I} g(\nub)}~
\left(\sum_{\nub \in \N_0^n} \frac{f(p^{\nu_1},\dots,p^{\nu_n})}{p^{\la \nub ,\cb\ra}}\right).$$
This ends the proof of point 3 and also the proof of Theorem \ref{main1}. \qed
\section{Proofs of Theorems \ref{main2} and \ref{main3}}
\subsection{Proof of Theorem \ref{main2}}
We will now explain how the combination of our Theorem \ref{main1} and La Bret\`eche's  multivariable Tauberian Theorem (i.e Theorems 1 and 2 of \cite{bretechecompo}) yields to our Theorem \ref{main2}. Our notations are different from  La Bret\`eche's notations. To simplify the exposition, we will first recall La Bret\`eche's Tauberian Theorem 1 and the part we use of his Tauberian Theorem 2 by using our notations:\\
{\bf Theorem A: (Theorem 1 of \cite{bretechecompo}):}\\
{\it Let $f: \N^n\rightarrow \R_+$ be a nonnegative function and $F$ the associated Dirichlet's series defined by
$$F(\s)=F(s_1,\dots,s_n)=\sum_{m_1,\dots,m_n\geq 1}\frac{f(m_1,\dots,m_n)}{m_1^{s_1}\dots  m_n^{s_n}}.$$
Denote by $\mathcal L\mathcal R_n^+(\C)$ the set of $\C-$linear forms from $\C^n$ to $\C$ that are nonnegative on $(\R_+)^n$.\\
We assume that there exists $\cb =(c_1,\dots,c_n)\in (\R_+)^n$ such that:
\be
\item $F(\s)$ converges absolutely for $\s\in \C^n$ such that $\Re(s_i)>c_i$ $\forall i=1,\dots,n$;
\item There exist a finite family  $\mathcal L =\left(\ell^{(i)}\right)_{1\leq i\leq q}$ of nonzero elements of 
$\mathcal L\mathcal R_n^+(\C)$, a finite family  $\left(h^{(i)}\right)_{1\leq i\leq q'}$ of elements of 
$\mathcal L\mathcal R_n^+(\C)$ and $\delta_1, \delta_2, \delta_3>0$ such that the function $H$ defined by 
$$H(\s)= F(\c+\s)~\prod_{i=1}^q \ell^{(i)} (\s)$$ has {\bf holomorphic} continuation to the domain
$$\mathcal D(\delta_1,\delta_3):=\{\s\in \C^n \mid \Re \left(\ell^{(i)} (\s)\right)>-\delta_1 ~\forall i=1,\dots,q {\mbox { and }}
\Re\left(h^{(i)} (\s)\right)>-\delta_3 ~\forall i=1,\dots,q'\}$$
and verifies the estimate: 
for $\eps, \eps'>0$ we have uniformly in $\s\in \mathcal D(\delta_1-\eps',\delta_3-\eps')$
$$H(\s)\ll \prod_{i=1}^q \left(|\Im \left(\ell^{(i)} (\s)\right)|+1)\right)^{1-\delta_2 \min \left(0, \Re \left(\ell^{(i)} (\s)\right)\right)}
\left(1+(|\Im (s_1)|+\dots+|\Im (s_n)|)^\eps\right).$$
\ee
Set $J=J(\cb)=\left\{j\in \{1,\dots, n\} \mid c_j=0\right\}$. Denote by $w=\# J$ the cardinality  of the set $J$ and by $j_1<\dots<j_w$ its elements in increasing order. 
Define the $w$ linear forms $\ell^{(q+i)}$ $(1\leq i\leq w)$ by 
$\ell^{q+i}(\s)=\e_{j_i}^*(\s)=s_{j_i}$.\\
Then, for any $\betab =(\beta_1,\dots, \beta_n) \in (0,\infty)^n$, there exist a polynomial $Q_{\betab}\in \R[X]$ of degree at most $q+ w-Rank \left\{ \ell^{(1)},\dots, \ell^{(q+w)}\right\}$ and $\theta >0$ such that 
$$\sum_{1\leq m_1 \leq x^{\beta_1}}\dots  \sum_{1\leq m_n \leq x^{\beta_n}} f(m_1,\dots,m_n)=x^{\langle \cb, \betab\rangle} Q_{\betab}(\log x) + O(x^{\langle \cb, \betab\rangle -\theta}) {\mbox { as }} x\rightarrow \infty.$$
}\\
{\bf Theorem B: (parts (ii) and (iv) from Theorem 2 of \cite{bretechecompo}):}\\
{\it Let $f: \N^n\rightarrow \R_+$ be a function satisfying assumptions of Theorem A. \\
Let $\betab =(\beta_1,\dots, \beta_n) \in (0,\infty)^n$. Set $\mathcal B = \sum_{i=1}^n \beta_i \e_i^* \in \mathcal L\mathcal R_n^+(\C)$.
\begin{itemize}
\item {\bf (ii)} If the Dirichlet's series $F$ satisfies the additional two assumptions:\\
(C1) There exists a function $G$ such $H(\s)= G\left(\ell^{(1)}(\s),\dots, \ell^{(q+w)}(\s)\right)$.\\
(C2) $\mathcal B \in Vect \left( \{\ell^{(k)}\mid k=1,\dots, q+w\}\right)$ and there is no subfamily $\mathcal L'$ of 
$\mathcal L_0:= \left(\ell^{(k)}\right)_{1\leq k \leq  q+w}$ 
such that  $\mathcal L' \neq \mathcal L_0$, 
$\mathcal B \in Vect(\mathcal L')$ and 
$\# \mathcal L' - Rank (\mathcal L') = \# \mathcal L_0 - Rank (\mathcal L_0)$.\par
Then, the polynomial $Q_{\betab}$ satisfies the relation 
$$Q_{\betab} (\log x)= H(\zerob) x^{-\langle \cb, \betab\rangle} \mathcal I_{\betab} (x) + O\left((\log x)^{\rho -1}\right),$$
where $\rho :=q+ w-Rank \left\{ \ell^{(1)},\dots, \ell^{(q+w)}\right\}$ and 
$$\mathcal I_{\betab}(x) :=\int_{\mathcal A_{\betab} (x)} \frac{dy_1\dots dy_q}{\prod_{i=1}^q y_i^{1-\ell^{(i)}(\cb)}},$$
with 
$$\mathcal A_{\betab} (x):=\{ \y \in [1,\infty)^q\mid \prod_{i=1}^q y_i^{\ell^{(i)}(\e_j)}\leq x^{\beta_j} ~~\forall j=1,\dots,n\}.
$$
\item {\bf (iv)} If $Rank \left\{ \ell^{(1)},\dots, \ell^{(q+w)}\right\}=n$, $H(\zerob)\neq 0$ and $\mathcal B\in con^*\left(\left\{ \ell^{(1)},\dots, \ell^{(q+w)}\right\}\right)$, then $deg(Q_{\betab})=\rho = q+w-n$.
\end{itemize}
}
{\bf Remark :} If assumptions of point (iv) hold, then assumptions of the point (ii) also clearly hold.\par
{\bf Proof of Theorem \ref{main2}:}\\
Let $f:\N^n\to \R_+$ be a multivariable multiplicative function. We assume that $f$ belongs to the class  
$\mathcal C (g,\kappa, \c, \delta)$ associated to the data $(g,\kappa, \c, \delta)$ (see definitions 1 and 2).
We assume also that the finite set
$$
I=I(\kappa, g):=\{\nub \in \N_0^n \mid \kappa(\nub)=1{\mbox{ and }} g(\nub)\neq 0\} \quad {\mbox {is nonempty.}}
$$
We denote by $\nub^1,\dots,\nub^r$ the elements of $I$ where $r=\# I$,\\
and  define  the finite sequence $q_k$ $(0\leq k\leq r)$ by
$$
q_0=0 \quad {\mbox{ and }} \quad q_k=\sum_{j=1}^k g(\nub^j) ~\forall k=1,\dots, r.
$$
We define the linear forms $\ell^{(i)}$ $(1\leq i\leq q_r)$ by 
$$\ell^{(i)}(\s)= \langle \nub^k, \s\rangle \quad {\mbox { if }} q_{k-1}<i\leq q_k {\mbox{ and }} 1\leq k \leq r.$$
We define also the set $J=J(\cb)=\left\{j\in \{1,\dots, n\} \mid c_j=0\right\}$. We 
denote by $w=\# J$ the cardinality  of the set $J$ and by $j_1<\dots<j_w$ its elements in increasing order. \\
We define also the $w$ linear forms $\ell^{(q+i)}$ $(1\leq i\leq w)$ by 
$$\ell^{q+i}(\s)=\e_{j_i}^*(\s)=s_{j_i} \quad (1\leq i\leq w).$$
By using notation of our Theorem \ref{main1} it's easy to see that the Dirichlet's series associated to $f$ is 
$$F(\s)=\sum_{m_1,\dots,m_n\geq 1}\frac{f(m_1,\dots,m_n)}{m_1^{s_1}\dots  m_n^{s_n}}=\mathcal M(f;\s)$$
and 
$$H(\s)= \left(\prod_{i=1}^{q_r} \ell^{(i)} (\s)\right)~F(\c+\s) =\left(\prod_{\nub \in I} \la \nub, \s\ra^{g(\nub)}\right)~\mathcal M (f;\c+\s)=\mathcal H(f,\c;\s).$$
Our Theorem \ref{main1} implies then that $F(\s)$ converges absolutely if $\Re(s_i)>c_i$ $\forall i=1,\dots,n$ and that 
there exists $\eps_0>0$ such that the function $\s\to H(\s)$
has holomorphic continuation to the domain $\{\s\in \C^n \mid \Re (s_i) > -\eps_0 ~\forall i=1,\dots,n\}$ and verifies in it the following estimate: for all $\eps >0$,
$$\mathcal H(f,\c;\s) \ll_\eps \prod_{\nub \in I} \left(|\la \nub, \s\ra|+1\right)^{g(\nub)\left(1-\frac{1}{2}\min\left(0, \Re (\la \nub, \s\ra)\right)\right)+\eps}.$$
For $i\in \{1,\dots,n\}$ set $h^{(i)}(\s)=s_i$ for all $\s=(s_1,\dots,s_n)\in \C^n$. Set also $\delta_1=\delta_3=\eps_0$, $q=q_r$ and $q'=n$. It follows then that  $\s\to H(\s)$ has 
{\bf holomorphic} continuation to the domain $$\mathcal D(\delta_1,\delta_3):=\{\s\in \C^n \mid \Re \left(\ell^{(i)} (\s)\right)>-\delta_1 ~\forall i=1,\dots,q {\mbox { and }}
\Re\left(h^{(i)} (\s)\right)>-\delta_3 ~\forall i=1,\dots,q'\}$$
and verifies the estimate: 
for $\eps, \eps'>0$ we have uniformly in $\s\in \mathcal D(\delta_1-\eps',\delta_3-\eps')$
$$H(\s)\ll \prod_{i=1}^q \left(|\Im \left(\ell^{(i)} (\s)\right)|+1)\right)^{1-\delta_2 \min \left(0, \Re \left(\ell^{(i)} (\s)\right)\right)}
\left(1+(|\Im (s_1)|+\dots+|\Im (s_n)|)^\eps\right),$$
where $\delta_2=1/2$.
Thus, all the assumptions of Theorem A above hold. By applying Theorem A with $\betab =\unb =(1,\dots,1)$, we deduce that there exist a polynomial $Q_{\unb}$ of degree at most 
$$\rho= q_r+ w-Rank \left\{ \ell^{(1)},\dots, \ell^{(q+w)}\right\} =\left(\sum_{\nub \in I} g(\nub)\right) +\#J -Rank\left(I\cup J\right)$$
and a positive constant $\eta >0$ such that
$$\mathcal N_\infty (f; x):=\sum_{\m=(m_1,\dots,m_n)\in \N^n\atop \|\m\|_\infty=\max_i m_i \leq x} f(m_1,\dots,m_n)
=x^{\|\c\|_1} Q_{\unb}(\ln x) + O\left(x^{\|\c\|_1-\eta}\right)\quad {\mbox{ as }} x\rightarrow \infty.$$
This ends the proof of the first part of our Theorem \ref{main2}.\par
Assume now  in addition that the two following assumptions hold:
\be
\item $Rank\left(I\cup J\right) =n$;
\item $\unb=(1,\dots,1)$ is in the interior of the cone generated by $I\cup J$; that is $ \unb \in con^*\left(I\cup J\right):=
\{\sum_{\nub \in I\cup J}\lambda_{\nub} \nub \mid \lambda_{\nub} \in (0,\infty)~\forall \nub \in I\cup J\}$,
\ee
By duality, we deduce that $Rank \left\{ \ell^{(1)},\dots, \ell^{(q+w)}\right\}=n$ and $\mathcal \unb^*\in con^*\left(\left\{ \ell^{(1)},\dots, \ell^{(q+w)}\right\}\right)$, 
Moreover since $f$ is nonegative, our Theorem \ref{main1} implies that 
$$H(\zerob)=\mathcal H(f,\c;\zerob)=\prod_p \left(1-\frac{1}{p}\right)^{\sum_{\nub \in I} g(\nub)} \left(\sum_{\nub \in \N_0^n}
\frac{f(p^{\nu_1},\dots,p^{\nu_n})}{p^{\la \nub, \c\ra}}\right)>0.$$
It follows that assumptions of point (iv) (and therefore assumptions of point (ii)) of Theorem B above hold. Theorem B implies then that 
$$deg(Q_{\unb})=\rho = \left(\sum_{\nub \in I} g(\nub)\right) +\#J -n$$
and 
\begin{equation}\label{Qconstant}
Q_{\unb} (\log x)= H(\zerob) x^{-\|\c\|_1} \mathcal I_{\unb} (x) + O\left((\log x)^{\rho -1}\right),
\end{equation}
where 
$$\mathcal I_{\unb}(x) =\int_{\mathcal A_{\unb} (x)} \frac{dy_1\dots dy_{q_r}}{\prod_{i=1}^{q_r} y_i^{1-\ell^{(i)}(\cb)}},$$
with 
$$\mathcal A_{\unb} (x):=\{ \y \in [1,\infty)^{q_r}\mid \prod_{i=1}^{q_r} y_i^{\ell^{(i)}(\e_j)}\leq x ~~\forall j=1,\dots,n\}.$$
By using notations of Definition \ref{Infxdef}, it's easy to see that 
$$\mathcal I_{\unb}(x) =\mathcal I_n(I,\u, \c; x) \quad {\mbox{ and }} \quad \mathcal A_{\unb} (x) =  \mathcal A (I,\u;x),$$
where $\u$ is the sequence $\u=\left(g(\nub)\right)_{\nub \in I}$. \\
Since the degree of the polynomial $Q_{\unb}$ is equal to $\rho=\left(\sum_{\nub \in I} g(\nub)\right) +\#J -n$, there exists a positive constant $C>0$ such that $Q_{\unb}(x)=C x^\rho +  O\left(x^{\rho -1}\right)$ as $x\rightarrow \infty$ and (\ref{Qconstant}) implies that 
$$ H(\zerob) x^{-\|\c\|_1} \mathcal I_{\unb} (x) =C (\log x)^\rho + O\left((\log x)^{\rho -1}\right).$$
It follows that 
$$C= H(\zerob)\lim_{x\rightarrow \infty}  x^{-\|\c\|_1} (\log x)^{-\rho}~\mathcal I_{\unb} (x)= H(\zerob)\lim_{x\rightarrow \infty}  x^{-\|\c\|_1} (\log x)^{-\rho}~\mathcal I_n(I,\u, \c; x).$$
We deduce that the main term of $\mathcal N_\infty (f; x)$ is given by
$$\mathcal N_\infty (f; x) = C_n(f)  K_n(f,\| \|_\infty) ~x^{\|\c\|_1} (\ln x)^{\rho}+ O\left((\ln x)^{\rho-1}\right) \quad {\mbox{as }} x\rightarrow \infty,$$
where $C_n(f):= H(\zerob)=\mathcal H(f,\c;\zerob)>0$ is defined by the Euler product (\ref{c0}) and
$$K_n(f,\| \|_\infty):=\lim_{x\rightarrow \infty} {\mathcal I}_n(I,\u,\c; x) ~x^{-\|\c\|_1} (\ln x)^{-\rho} >0.$$
This ends the proof of Theorem \ref{main2}.\qed
\subsection{Proof of Theorem \ref{main3}}
Let $f:\N^n\to \R$ be a multivariable multiplicative function. We assume that $f$ belongs to the class  
$\mathcal C (g,\kappa, \c, \delta)$ associated to the data $(g,\kappa, \c, \delta)$ (see definitions 1 and 2).\\
We assume also that the finite set
$$
I=I(\kappa, g):=\{\nub \in \N_0^n \mid \kappa(\nub)=1{\mbox{ and }} g(\nub)\neq 0\} \quad {\mbox {is nonempty.}}
$$
We define the set $I_{\cb}:=\{\frac{1}{\langle\nub, \c\rangle} \nub\mid \nub \in I\}$ and the sequence $\u:=\left(u(\betab)\right)_{\betab \in I_{\cb}}$ where 
$$u(\betab)= \sum_{\nub \in I;~ \frac{1}{\langle\nub, \c\rangle} \nub =\betab}   g(\nub) \quad {\mbox{ for all }}\quad  \betab \in I_{\cb}.$$
We Define also the pair $\T_{\cb}:=(I_{\cb}, \u)$.\\
Theorem \ref{main1} implies that
$$\s\to \mathcal M (f;\s):=\sum_{m_1\geq 1,\dots,m_n\geq 1}\frac{f(m_1,\dots,m_n)}{m_1^{s_1}\dots m_n^{s_n}}$$
converges absolutely in the domain $\{\s\in \C^n \mid \Re (s_i) > c_i ~\forall i=1,\dots,n\}$;
and that there exists $\eps_0>0$ such that the function 
\begin{eqnarray}\label{Hfacteur}
\s\to H(f;\T_{\cb}; \s)&:=&\left(\prod_{\betab \in I_{\cb}} \la \betab, \s\ra^{u(\betab)}\right)~\mathcal M (f;\c+\s)\nonumber\\
&=&\left(\prod_{\nub \in I} \la \nub, \cb\ra^{-g(\nub)}\right) \left(\prod_{\nub \in I} \la \nub, \s\ra^{g(\nub)}\right)~
\mathcal M (f;\c+\s)\nonumber\\
&=& \left(\prod_{\nub \in I} \la \nub, \cb\ra^{-g(\nub)}\right) ~\mathcal H(f,\c;\s)
\end{eqnarray}
has {\it holomorphic} continuation to the domain $\{\s\in \C^n \mid \Re (s_i) > -\eps_0 ~\forall i=1,\dots,n\}$ and verifies in it the following estimate: for all $\eps >0$,
$$H(f,\T_{\cb};\s) \ll_\eps \prod_{\nub \in I} \left(|\la \nub, \s\ra|+1\right)^{g(\nub)\left(1-\frac{1}{2}\min\left(0, \Re (\la \nub, \s\ra)\right)\right)+\eps};$$
We deduce that $f$ is of {\it finite type} with $\T_{\cb}:=(I_{\cb}, \u)$ as a {\it regularizing pair}  (see Definition 2 of \cite{Toricmanin}). It follows then from Corollary 2 of \cite{Toricmanin} that 
there exist a polynomial $Q$ of degree at most 
$$\rho:=\left(\sum_{\betab \in I_{\cb}} u(\betab)\right)  -Rank (I_{\cb})=\left(\sum_{\nub \in I} g(\nub)\right) - Rank (I)$$
 and a positive constant $\mu >0$ such that
$$\mathcal N_d (f; x):=\sum_{\m=(m_1,\dots,m_n)\in \N^n\atop \|\m\|_d =\sqrt[d]{m_1^d+\dots+m_n^d}\leq x}
f(m_1,\dots,m_n)
=x^{\|\c\|_1} Q(\ln x) + O\left(x^{\|\c\|_1-\mu}\right)\quad {\mbox{ as }} x\rightarrow \infty.$$
This ends the proof of part 1 of Theorem \ref{main3}.\par
Assume now  in addition that $Rank ( I )=n$ and $ \unb \in con^*\left(I\right)$. It follows that 
\be
\item $Rank ( I_{\cb} )=n$ and it's clear then that there exists a function holomorphic in a tubular neighborhood of $\zerob$ such that $H(f,\T_{\cb};\s)=K\left( (\langle \betab, \s\rangle)_{\betab \in I_{\cb}}\right)$;
\item $ \unb \in con^*\left(I_{\cb}\right)$. 
 \ee 
 Therefore, the additional assumptions 1 and 2 of Theorem 3 of \cite{Toricmanin} are satisfied and the second part of Corollary 2 of \cite{Toricmanin} implies then that 
 $$\mathcal N_d (f; x) = C_0(f, P_d) ~x^{\|\c\|_1} (\ln x)^{\rho}+ O\left((\ln x)^{\rho-1}\right) \quad {\mbox{as }} x\rightarrow \infty,$$
 where $\ds C_0(f, P_d):=\frac{H(f,\T_{\cb};\zerob) d^{\rho+1} ~A_0(\T_{\cb}, P_d)}{\|\cb\|_1~ \rho !}$
 and $A_0(\T_{\cb}, P_d)>0$ is the mixed volume constant (see \S  2.3.3  of \cite{Toricmanin}) associated to the pair 
 $\mathcal T_{\cb} :=(I_{\cb}, \u)$ and the polynomial $P_d= X_1^{d}+\dots + X_n^{d}$.\\
 Combining (\ref{Hfacteur}) and the expression of $\mathcal H(f,\c;\zerob)$ given by theorem \ref{main1} implies that 
 $$H(f,\T_{\cb};\zerob)=\left(\prod_{\nub \in I} \la \nub, \cb\ra^{-g(\nub)}\right) ~C_n(f),$$
 where $C_n(f):= \mathcal H(f,\c;\zerob)>0$ is defined by the Euler product (\ref{c0}).
 Moreover, if we set 
$$ K_n(f,\| \|_d):=\left(\prod_{\nub \in I} \la \nub, \cb\ra^{-g(\nub)}\right)~\frac{d^{\rho+1}~
 A_0(\T_{\cb} ,P_d)}{\|\c\|_1 ~ \rho!} > 0,$$
 then the the constant  $C_0(f, P_d)$ is positive and is given by 
 $$C_0(f, P_d)=C_n(f) ~ K_n(f,\| \|_d)>0.$$
In particular, the degree of the polynomial $Q$ is equal to $\rho=\left(\sum_{\nub \in I} g(\nub)\right) -n$.
This ends the proof of Theorem \ref{main3}.\qed
\section{Proof of Corollary \ref{corollary1}}
Define the function $g_1:\N_0^n\rightarrow \N_0$ by
$$g_1(\nub)= \left\{
    \begin{array}{ll}
        1 \qquad \qquad  {\mbox{ if }}  \exists  i\neq j\in \{1,\dots, n\} {\mbox { such that }}\nu_i=\nu_j=\|\nub\|_\infty;\\
       \|\nub\|_\infty -\max\left(\{\nu_i\mid i=1,\dots,n\}\setminus \{\|\nub\|_\infty\}\right)+1
        \quad \mbox{ otherwise,}
    \end{array}
\right.
$$
where $ \|\nub\|_\infty =\max_{i=1,\dots,n} \nu_i$.\\
We will first prove the following needed lemma.
\begin{lemma}\label{neededlemma}
We have
$$c_n(p^{\nu_1},\dots, p^{\nu_n})= g_1(\nub)~ p^{\|\nub\|_1-\|\nub\|_\infty} +O\left((1+\|\nub\|_1)
 p^{\|\nub\|_1-\|\nub\|_\infty -1} \right)$$
 uniformly in $\nub=(\nu_1,\dots,\nu_n)\in \N_0^n$ and $p$ prime number.
\end{lemma}
{\bf Proof of Lemma \ref{neededlemma}:}\\
In the proof of this lemma we will use the notations: $a\wedge b= \min(a,b)$ and $a\vee b=\max(a,b)$.\\
First we recall the following formula proved by T\'oth in \cite{T2012}:
\begin{equation}\label{TF}
c_n(p^{\nu_1},\ldots,p^{\nu_n})=
\underset{\substack{0\leq\ell_i\leq v_i\\i\in\llbracket 1,n\rrbracket}}{\displaystyle\sum}
\dfrac{\varphi(p^{\ell_1})\cdots \varphi(p^{\ell_n})}{\varphi\left(p^{\max\{\ell_1,\ldots,\ell_n\}}\right)},
\end{equation}
where $\varphi$ is the Euler's  totient function.\\
If $n=1$, then $c_1(p^{\nu_1})=1+\nu_1$ and the lemma holds. \\
Let $n\geq 2$. Set $k=n-1\in \N$. \\
Let $p$ be a prime number and $\nub =(\nu_1,\dots,\nu_{n}) \in \N_0^{n}$. Without loss of generality we can assume that
$$\nu_1\leq \nu_2 \leq \dots \leq \nu_{n}.$$
It follows that
\begin{align*}
c_{n}(p^{\nu_1},\ldots,p^{\nu_{n}})&=
\underset{\substack{0\leq\ell_i\leq \nu_i\\i\in\llbracket 1,k+1\rrbracket\\\ell_1\vee\cdots\vee\ell_k\leq
\ell_{k+1}}}{\displaystyle\sum}
\varphi(p^{\ell_1})\cdots\varphi(p^{\ell_{k}}) +
\underset{\substack{0\leq\ell_i\leq \nu_i\\i\in\llbracket 1,k+1\rrbracket\\\ell_1\vee\cdots\vee\ell_k>\ell_{k+1}}}{\displaystyle\sum}
\varphi(p^{\ell_{k+1}})
\dfrac{\varphi(p^{\ell_1})\cdots\varphi(p^{\ell_k})}{\varphi\left(p^{\max\{\ell_1,\ldots,\ell_k\}}\right)}\\
&={\displaystyle\sum\limits_{\ell=0}^{\nu_{k+1}}}\quad
p^{\nu_1\wedge\ell+\cdots+\nu_k\wedge\ell}
+\displaystyle\sum\limits_{\substack{0\leq\ell_i\leq \nu_i\\i\in\llbracket 1,k\rrbracket\\\ell_1\vee\cdots\vee\ell_k\geq1}}
\dfrac{\varphi(p^{\ell_1})\cdots\varphi(p^{\ell_k})}{\varphi\left(p^{\max\{\ell_1,\ldots,\ell_k\}}\right)}
\displaystyle\sum\limits_{\ell_{k+1}=0}^{\ell_1\vee\cdots\vee\ell_k-1}
\varphi(p^{\ell_{k+1}})\\
&=(\nu_{k+1}-\nu_k+1) p^{\nu_1+\cdots+\nu_k}
+{\displaystyle\sum\limits_{\ell=0}^{\nu_{k}-1}}\quad
p^{\nu_1\wedge\ell+\cdots+\nu_k\wedge\ell}\\
&+
\displaystyle\sum\limits_{\substack{0\leq\ell_i\leq \nu_i\\i\in\llbracket 1,k\rrbracket\\\ell_1\vee\cdots\vee\ell_k\geq1}}
\dfrac{\varphi(p^{\ell_1})\cdots\varphi(p^{\ell_k})}{\varphi\left(p^{\max\{\ell_1,\ldots,\ell_k\}}\right)}p^{\ell_1\vee\cdots\vee\ell_k-1}\\
&=(\nu_{k+1}-\nu_k+1)p^{\nu_1+\cdots+\nu_k}
+{\displaystyle\sum\limits_{\ell=0}^{\nu_{k}-1}}\quad
p^{\nu_1\wedge\ell+\cdots+\nu_k\wedge\ell}
+
\displaystyle\sum\limits_{\substack{0\leq\ell_i\leq \nu_i\\i\in\llbracket 1,k\rrbracket}}
\dfrac{\varphi(p^{\ell_1})\cdots\varphi(p^{\ell_k})}{p-1}-\dfrac{1}{p-1}
\end{align*}
\begin{align*}
\hspace{2.5cm}&=(\nu_{k+1}-\nu_k+1)p^{\nu_1+\cdots+\nu_k}
+{\displaystyle\sum\limits_{\ell=0}^{\nu_{k}-1}}\quad
p^{\nu_1\wedge\ell+\cdots+\nu_k\wedge\ell}
+
\dfrac{p^{\nu_1+\cdots+\nu_k}-1}{p-1}\\
&=(\nu_{k+1}-\nu_k+1)p^{\nu_1+\cdots+\nu_k}
+{\displaystyle\sum\limits_{\ell=0}^{\nu_{k}-1}}\quad
p^{\nu_1\wedge\ell+\cdots+\nu_{k-1}\wedge\ell} p^{\ell}+\sum_{\ell=0}^{\nu_1+\cdots+\nu_k-1}p^{\ell}
\end{align*}
Thus, for $\nub =(\nu_1,\dots,\nu_{n}) \in \N_0^{n}$ such that $\nu_1\leq \nu_2 \leq \dots \leq \nu_{n}$, we have
\begin{equation}\label{cnpnuformulatouse}
c_{n}(p^{\nu_1},\ldots,p^{\nu_{n}}) =
(\nu_{n}-\nu_{n-1}+1)p^{\nu_1+\cdots+\nu_{n-1}}
+{\displaystyle\sum\limits_{\ell=0}^{\nu_{n-1}-1}}\quad
p^{\nu_1\wedge\ell+\cdots+\nu_{n-2}\wedge\ell} p^{\ell}+\sum_{\ell=0}^{\nu_1+\cdots+\nu_{n-1}-1}p^{\ell}
\end{equation}
We deduce that
\begin{eqnarray*}
&& 0\leq c_{n}(p^{\nu_1},\ldots,p^{\nu_{n}})-(\nu_{n}-\nu_{n-1}+1)
p^{\|\nub\|_1-\|\nub\|_{\infty}} \\
&\leq&  p^{\nu_1+\cdots+\nu_{k-1}} {\displaystyle\sum\limits_{\ell=0}^{\nu_{k}-1}}\quad
 p^{\ell} + (\nu_1+\dots+\nu_k) p^{\nu_1+\dots+\nu_k -1}\\
 &\leq&  (\nu_1+\dots+\nu_{k-1} +2\nu_k) p^{\nu_1+\dots+\nu_k -1}
 \leq 2\|\nub\|_1 ~ p^{\|\nub\|_1-\|\nub\|_{\infty}-1}.
\end{eqnarray*}
This ends the proof of Lemma \ref{neededlemma}. \qed

We will now use Lemma \ref{neededlemma} to prove Corollary \ref{corollary1}. \\
It's clear that  $c_n: (m_1,\dots,m_n)\mapsto c_n(m_1,\dots,m_n)$ is a multiplicative function. Moreover, Lemma \ref{neededlemma} implies that $c_n$ belongs to the class $\mathcal C (g,\kappa, \c, \delta)$ (see definition \ref{class-def}), where
$g=g_1$, $\c=\unb=(1,\dots,1)$, $\delta=1$ and $\kappa$ is the function defined by $\kappa (\nub) =\max_{i=1,\dots,n} \nu_i$
$\forall \nub \in \N_0^n$.
Furthermore, if we denote by $(\e_1,\dots, \e_n)$ the canonical basis of $\R^n$, then
$$I=I(\kappa, g):=\{\nub \in \N_0^n \mid \kappa(\nub)=1 {\mbox{ and }} g(\nub)\neq 0\}
=\{0,1\}^n\setminus\{\zerob\}.$$
Since $J=\{\e_i\mid c_i=0\}=\emptyset$ and $\e_1,\dots, \e_n \in I=I\cup J$, it follows that the two assumptions
$Rank ( I\cup J )=n$ and $ \unb \in con^*\left(I\cup J\right)$ hold.
Set
$$\rho:=\left(\sum_{\nub \in I} g(\nub)\right) +\#J -n=\left(\sum_{\nub \in I} g(\nub)\right) -n.$$
Since $g(\e_i)=2$ $\forall i=1,\dots,n$ and
$g(\nub) =1$ $\forall \nub \in I\setminus \{e_1,\dots, e_n\}$, we have
$$\rho= 2n+ (\#I -n) -n=\# I =2^n -1.$$
Theorem \ref{main2} implies then that there exist a polynomial $Q_1$ of degree $\rho$ and a positive constant $\mu_1 >0$ such that
\begin{eqnarray*}
G_n(x) &:=&\sum_{1\leq m_1,\dots, m_n \leq x}c_n(m_1,\dots,m_n)
=x^{n} Q_1(\ln x) + O\left(x^{n-\mu_1}\right)\quad {\mbox{ as }} x\rightarrow \infty,\\
&=& C_n(c_n)  K_n(c_n,\| \|_\infty)~x^n (\ln x)^{2^n-1} +O\left( x^n (\ln x)^{2^n-2}\right) \quad {\mbox{ as }} \quad x\rightarrow \infty,
\end{eqnarray*}
where
$$C_n(c_n):=  \mathcal H(c_n,\c;\zerob)=\prod_p \left(1-\frac{1}{p}\right)^{2^n+n-1} \left(\sum_{\nub \in \N_0^n}
\frac{c_n(p^{\nu_1},\dots,p^{\nu_n})}{p^{\|\nub\|_1}}\right) >0$$
 and
 $$K_n(c_n,\| \|_\infty):=\lim_{x\rightarrow \infty} {\mathcal I}_n(I,\u; x) ~x^{-n} (\ln x)^{-2^n+1} >0,\qquad {\mbox{ where}}
 $$
$\ds \mathcal I_n(I,\u,\c; x)$ is the integral (see definition \ref{Infxdef}) associated to the set $\ds I=\{0,1\}^n\setminus \{\zerob\}$ and to the sequence
$\ds \u=\left(u(\nub)\right)_{\nub \in I}$ defined by
$u(\e_i)=2$ $\forall i=1,\dots,n$ and $u(\nub) =1$ $\forall \nub \in I\setminus \{e_1,\dots, e_n\}$ and to the vector $\c=\unb$.
This ends the proof of corollary \ref{corollary1}. \qed

\section{Proof of Corollaries \ref{corollary2}, \ref{corollary3}, \ref{corollary4},  \ref{corollary5} and  \ref{corollary6}}
\subsection{Proof of Corollary \ref{corollary2}}
Let $s_n:\N^n\rightarrow \R_+$ be the function defined by
$$s_n(m_1,\dots,m_n) =  \frac{1}{\lcm(m_1,\dots,m_n)}\quad \forall (m_1,\dots,m_n) \in \N^n.$$
It clear that the function $s_n$ is multiplicative and that for $\nub=(\nu_1,\dots,\nu_n)\in \N_0^n$ and $p$ prime number, we have
$\ds s_n(p^{\nu_1},\dots, p^{\nu_n})= p^{-\max_{i=1,\dots,n} \nu_i}.$
Thus,  $s_n$ belongs to the class $\mathcal C (g,\kappa, \c, \delta)$ (see definition \ref{class-def}), where
$g\equiv 1$, $\c=\zerob=(0,\dots,0)$, $\delta=1$ and $\kappa$ is the function defined by $\kappa (\nub) =\max_{i=1,\dots,n} \nu_i$ $\forall \nub \in \N_0^n$.

Moreover, we have
$$I=I(\kappa, g):=\{\nub \in \N_0^n \mid \kappa(\nub)=1 {\mbox{ and }} g(\nub)\neq 0\}
=\{0,1\}^n\setminus \{\zerob\}$$
and $J=\{\e_i\mid c_i=0\}=\{\e_1,\dots,\e_n\}$. It follows that the two assumptions
$Rank ( I\cup J )=n$ and $ \unb \in con^*\left(I\cup J\right)$ hold.
Moreover, $\rho:=\left(\sum_{\nub \in I} g(\nub)\right) +\#J -n=\left(\sum_{\nub \in I} g(\nub)\right) = 2^n-1 $.\\
Theorem \ref{main2} implies then that there exist a polynomial $Q_2$ of degree $\rho=2^n-1$ and a positive constant $\mu_2 >0$ such that
\begin{eqnarray*}
S_n(x)&:=&\sum_{1\leq m_1,\dots, m_n \leq x} \frac{1}{\lcm(m_1,\dots,m_n)} = Q_2(\ln x) + O(x^{-\mu_2}) \quad {\mbox{ as }} \quad x\rightarrow \infty,\\
&=& C_n(s_n)  K_n(s_n,\| \|_\infty)~(\ln x)^{2^n-1} +O\left( (\ln x)^{2^n-2}\right) \quad {\mbox{ as }} \quad x\rightarrow \infty,
\end{eqnarray*}
where
\begin{eqnarray*}
C_n(s_n)&:=&\mathcal H(c_n,\c;\zerob)=\prod_p \left(1-\frac{1}{p}\right)^{2^n-1} \left(\sum_{\nub \in \N_0^n}
\frac{1}{p^{\|\nub\|_\infty}}\right) \\
&=&\prod_p \left(1-\frac{1}{p}\right)^{2^n-1} \left(\sum_{k=0}^\infty
\frac{(k+1)^n -k^n}{p^{k}}\right)>0,
\end{eqnarray*}
 $${\mbox{and}}\quad K_n(s_n,\| \|_\infty):=\lim_{x\rightarrow \infty} {\mathcal I}_n(I,\u,\c; x) ~(\ln x)^{-2^n+1} >0,\qquad {\mbox{ where}}$$
 $\ds \mathcal I_n(I,\u, \c; x)$ is the integral (see definition \ref{Infxdef}) associated to $\ds I=\{0,1\}^n\setminus \{\zerob\}$,  to the sequence
$\ds \u=\left(u(\nub)\right)_{\nub \in I}$ defined by
$u(\nub) =1$ $\forall \nub \in I$
and to the vector $\c=\zerob$. This ends the proof of corollary \ref{corollary2}. \qed

\subsection{Proof of Corollary \ref{corollary3}}
Let $n \in \N\setminus\{1\}$. Let $u_n:\N^n\rightarrow \R_+$ be the function defined by
$$u_n(m_1,\dots,m_n) =  \frac{1}{\lcm(m_1,\dots,m_n)} {\mbox{ if }} \gcd(m_1,\dots,m_n)=1 {\mbox { and }} u_n(m_1,\dots,m_n) =0 {\mbox{ otherwise}}.$$
It is clear that the function $u_n$ is multiplicative and that for $\nub=(\nu_1,\dots,\nu_n)\in \N_0^n$ and $p$ prime number, we have
$$u_n(p^{\nu_1},\dots, p^{\nu_n})= p^{-\max_{i=1,\dots,n} \nu_i} {\mbox{ if }} \min_{i=1,\dots,n} \nu_i =0 {\mbox { and }} u_n(p^{\nu_1},\dots, p^{\nu_n})=0 {\mbox{ otherwise}}.$$
Thus,  $u_n$ belongs to the class $\mathcal C (g,\kappa, \c, \delta)$ (see definition \ref{class-def}), where
$\c=\zerob=(0,\dots,0)$, $\delta=1$, $\kappa$ is the function defined by $\kappa (\nub) =\max_{i=1,\dots,n} \nu_i$ $\forall \nub \in \N_0^n$ and $g$ is the function defined by
$$g(\nub)=1 {\mbox{ if }} \min_{i=1,\dots,n} \nu_i =0 {\mbox { and }} g(\nub)=0 {\mbox { otherwise}}.$$
Thus, we have
$\ds I=I(\kappa, g):=\{\nub \in \N_0^n \mid \kappa(\nub)=1 {\mbox{ and }} g(\nub)\neq 0\}
=\{0,1\}^n\setminus\{\zerob, \unb\}$ \\
and $J=\{\e_i\mid c_i=0\}=\{\e_1,\dots,\e_n\}$. It follows that the two assumptions
$Rank ( I\cup J )=n$ and $ \unb \in con^*\left(I\cup J\right)$ hold.
Moreover, $\rho:=\left(\sum_{\nub \in I} g(\nub)\right) +\#J -n=\left(\sum_{\nub \in I} g(\nub)\right) = 2^n-2 $.\\
Theorem \ref{main2} implies then that there exist a polynomial $Q_3$ of degree $2^n-2$ ans $\mu_3 >0$ such that
\begin{eqnarray*}
U_n(x)&:=&\sum_{1\leq m_1,\dots, m_n \leq x\atop \gcd(m_1,\dots,m_n)=1} \frac{1}{\lcm(m_1,\dots,m_n)} = Q_3(\ln x) +
 O(x^{-\mu_3}) \quad {\mbox{ as }} \quad x\rightarrow \infty\\
 &=& C_n(u_n)  K_n(u_n,\| \|_\infty)~(\ln x)^{2^n-2} +O\left( (\ln x)^{2^n-3}\right) \quad {\mbox{ as }} \quad x\rightarrow \infty,
\end{eqnarray*}
where
\begin{eqnarray*}
C_n(u_n)&:=& \mathcal H(c_n,\c;\zerob)=\prod_p \left(1-\frac{1}{p}\right)^{2^n-2} \left(\sum_{\nub \in \N_0^n\atop
 \min_{i=1,\dots,n} \nu_i =0}
\frac{1}{p^{\|\nub\|_\infty}}\right) \\
&=&\prod_p \left(1-\frac{1}{p}\right)^{2^n-2} \left(1+\sum_{k=1}^\infty
\frac{(k+1)^n +(k-1)^n-2k^n}{p^{k}}\right)\\
&=&\prod_p \left(1-\frac{1}{p}\right)^{2^n-1} \left(\sum_{k=0}^\infty
\frac{(k+1)^n -k^n}{p^{k}}\right)>0,
\end{eqnarray*}
 $${\mbox{and}} \quad K_n(u_n,\| \|_\infty):=\lim_{x\rightarrow \infty} {\mathcal I}_n(I,\u,\c; x) ~(\ln x)^{-2^n+2} >0,\qquad {\mbox{ where}}
 $$
$\ds \mathcal I_n(I,\u, \c; x)$ is the integral (see definition \ref{Infxdef}) associated to $\ds I=\{0,1\}^n\setminus \{\zerob, \unb\}$,  to the sequence
$\ds \u=\left(u(\nub)\right)_{\nub \in I}$ defined by
$u(\nub) =1$ $\forall \nub \in I$
and to the vector $\c=\zerob$.
This ends the proof of corollary \ref{corollary3}. \qed
\subsection{Proof of Corollary \ref{corollary4}}
Let $v_n:\N^n\rightarrow \R_+$ be the function defined by
$$v_n(m_1,\dots,m_n) =  \frac{m_1\dots m_n}{\lcm(m_1,\dots,m_n)}\quad \forall (m_1,\dots,m_n) \in \N^n.$$
It is clear that the function $f$ is multiplicative and that for $\nub=(\nu_1,\dots,\nu_n)\in \N_0^n$ and $p$ prime number, we have
$$v_n(p^{\nu_1},\dots, p^{\nu_n})= p^{\|\nub\|_1-\max_{i=1,\dots,n} \nu_i}.$$
Thus,  $v_n$ belongs to the class $\mathcal C (g,\kappa, \c, \delta)$ (see definition \ref{class-def}), where
$g\equiv 1$, $\c=\unb=(1,\dots,1)$, $\delta=1$ and $\kappa$ is the function defined by $\kappa (\nub) =\max_{i=1,\dots,n} \nu_i$ $\forall \nub \in \N_0^n$.

Moreover, we have
$$I=I(\kappa, g):=\{\nub \in \N_0^n \mid \kappa(\nub)=1 {\mbox{ and }} g(\nub)\neq 0\}
=\{0,1\}^n\setminus\{\zerob \}$$
and $J=\{\e_i\mid c_i=0\}=\emptyset$. It follows that the two assumptions
$Rank ( I\cup J )=n$ and $ \unb \in con^*\left(I\cup J\right)$ hold.
Moreover, $\rho:=\left(\sum_{\nub \in I} g(\nub)\right) +\#J -n= 2^n-1 -n$.\\
Theorem \ref{main2} implies then that there exist a polynomial $Q_4$ of degree $\rho=2^n-n-1$ and a positive constant $\mu_4 >0$ such that
\begin{eqnarray*}
V_n(x)&:=&\sum_{1\leq m_1,\dots, m_n \leq x} \frac{m_1\dots,m_n}{\lcm(m_1,\dots,m_n)} = x^n Q_4(\ln x) + O(x^{n-\mu_4}) \quad {\mbox{ as }} \quad x\rightarrow \infty\\
&=& C_n(v_n)  K_n(v_n,\| \|_\infty)~x^n~(\ln x)^{2^n-n-1} +O\left( x^n~(\ln x)^{2^n-n-2}\right) \quad {\mbox{ as }} \quad x\rightarrow \infty,
\end{eqnarray*}
where
\begin{eqnarray*}
C_n(v_n)&:=&\mathcal H(c_n,\c;\zerob)=\prod_p \left(1-\frac{1}{p}\right)^{2^n-1} \left(\sum_{\nub \in \N_0^n}
\frac{1}{p^{\|\nub\|_\infty}}\right)\\
 &=&\prod_p \left(1-\frac{1}{p}\right)^{2^n-1} \left(\sum_{k=0}^\infty
\frac{(k+1)^n -k^n}{p^{k}}\right)>0,
\end{eqnarray*}
 $${\mbox{and}}\quad K_n(v_n,\| \|_\infty):=\lim_{x\rightarrow \infty} {\mathcal I}_n(I,\u,\c; x) ~x^{-n} (\ln x)^{-2^n+n+1} >0,
 \qquad {\mbox{ where}}
 $$
$\ds \mathcal I_n(I,\u, \c; x)$ is the integral (see definition \ref{Infxdef}) associated to $\ds I=\{0,1\}^n\setminus \{\zerob\}$,  to the sequence
$\ds \u=\left(u(\nub)\right)_{\nub \in I}$ defined by
$u(\nub) =1$ $\forall \nub \in I$
and to the vector $\c=\unb$.
This ends the proof of corollary \ref{corollary4}. \qed
\subsection{Proof of Corollaries \ref{corollary5} and \ref{corollary6}}
Proof of corollary \ref{corollary5} (resp. corollary \ref{corollary6}) is similar to the proof of corollary \ref{corollary1} (resp. corollary \ref{corollary4}) by using Theorem \ref{main3} instead of Theorem \ref{main2} and the identity 
$ \prod_{\nub \in \{0,1\}^n,~\|\nub\|_1\geq 2} \|\nub\|_1 =\prod_{k=2}^n k^{{n\choose k}}$.

\section{Explicit computations of the constants $C_n(.)$ and $K_n(.,.)$ in dimensions  $n=2$ and $n=3$}
We will use the software WX Maxima to compute some iterated integrals below.
\subsection{Computation of $C_2(c_2)$ in corollaries \ref{corollary1} and \ref{corollary5}}
The identity (\ref{cnpnuformulatouse}) implies that  for $\nub \in \N_0^2$ such that $0\leq \nu_1\leq \nu_2$ we have
$$c_2(p^{\nu_1},p^{\nu_2})=(\nu_2-\nu_1-1)p^{\nu_1}+2\frac{p^{\nu_1+1}-1}{p-1}.$$
   We deduce by  symmetry  that
$$\begin{array}{cl}
C_2(c_2)&=\displaystyle\prod_{p}\left(1-\frac{1}{p}\right)^5\left[2\sum_{\nu_2\geq0}\sum_{\nu_1=0}^{\nu_2}\frac{c_2(p^{\nu_1},p^{\nu_2})}{p^{\nu_1+\nu_2}}-\sum_{\nu_1\geq0}\frac{c_2(p^{\nu_1},p^{\nu_1})}{p^{2 \nu_1}}\right]\\
&=\displaystyle\prod_p \left(\frac{p-1}{p}\right)^5\left[2\frac{p^2(p^2+p+2)}{(p-1)^3(p+1)}-\frac{p(p^2+1)}{(p-1)^3}\right] 
=\displaystyle\prod_p \left(1-\frac{1}{p^2}\right)^2=\frac{1}{\zeta^2(2)}=\frac{36}{\pi^4}.
\end{array}$$
\subsection{Computation of $C_n(.)$ $(n=2,3)$ in corollaries \ref{corollary2}, \ref{corollary3}, \ref{corollary4} and \ref{corollary6}}
Constants $C_n(.)$ in corollaries \ref{corollary2}, \ref{corollary3}, \ref{corollary4} and \ref{corollary6} are equal. We will denote them by $C_n$ in this subsection. \\
$\bullet$ In dimension $n=2$, we have
%\begin{align*}
%C_2=\prod_{p}\left(1-\frac{1}{p}\right)^3\left(\sum_{k\geq 0}\frac{2k+1}{p^k}\right)
%&=\prod\limits_{p}\frac{(p-1)^3}{p^3}\left(\frac{2p}{(p-1)^2}+\frac{p}{p-1}\right)\\
%&=\prod_p\frac{p^2-1}{p^2}=\zeta(2)^{-1}=\frac{6}{\pi^2}
%\end{align*}
$$C_2=\prod_{p}\left(1-\frac{1}{p}\right)^3\left(\sum_{k\geq 0}\frac{2k+1}{p^k}\right)
=\prod\limits_{p}\frac{(p-1)^3}{p^3}\left(\frac{2p}{(p-1)^2}+\frac{p}{p-1}\right)
=\prod_p\frac{p^2-1}{p^2}=\zeta(2)^{-1}=\frac{6}{\pi^2}.$$
$\bullet$ In dimension $n=3$, we have
\begin{align*}
C_3&=\prod_{p}\left(1-\frac{1}{p}\right)^7\left(\sum_{k\geq 0}\frac{3k^2+3k+1}{p^k}\right)
=\prod\limits_{p}\frac{(p-1)^7}{p^7}\left(\frac{3p(p+1)}{(p-1)^3}+\frac{3p}{(p-1)^2}+\frac{p}{p-1}\right)\\
&=\prod_p\frac{(p-1)^4(p^2+4p+1)}{p^6}
=\prod_p\left(1-\frac{9}{p^2}+\frac{16}{p^3}-\frac{9}{p^4}+\frac{1}{p^6}\right)
\end{align*}
\subsection{Computation of $K_n(c_n,\| \|_\infty) $ $(n=2,3)$ in corollary \ref{corollary1}}

$\bullet$ In dimension $n=2$, we have:
$I=\{\e_1,\e_2,\e_1+\e_2\}$, $\u= (2,2,1)$ and $\c=\unb$. It follows from definition \ref{Infxdef} that
\begin{align*}
{\mathcal I}_2(I,\u,\c; x)&=\int\limits_{\substack{\boldsymbol{y}\in[1,+\infty[^5\\
y_1y_2y_5\leq x\\y_3y_4y_5\leq x}}y_5\,d\boldsymbol{y}=
\int\limits_{y_5\in[1,x[}y_5\Bigl(\int\limits_{\substack{y_1,y_2\in[1,+\infty[\\y_1y_2\leq x/y_5}}\,dy_1\,dy_2\Bigr)^2\,dy_5\\
&=\int\limits_{y_5\in[1,x[} y_5\left(\frac{x}{y_5}\ln\left(\frac{x}{y_5}\right)-\frac{x}{y_5}+1\right)^2\,dy_5\\
&=\frac{1}{3}x^2\ln^3(x)-x^2\ln^2(x)+x^2\ln(x)-2x\ln(x)+\frac{x^2}{2}-\frac{1}{2},
\end{align*}
which implies that
$$K_2(c_2,\| \|_\infty):=\lim_{x\rightarrow \infty} {\mathcal I}_2(I,\u,\c; x) ~x^{-2} (\ln x)^{-3}=\frac{1}{3}.$$

$\bullet$ In dimension $n=3$: We have $I=\{\e_1,\e_2,\e_3,\e_1+\e_2,\e_1+\e_3,\e_2+\e_3,\e_1+\e_2+\e_3\}$,
$\u=(2,2,2,1,1,1,1)$ and $\c=\unb$. It follows from definition \ref{Infxdef} that
$$
{\mathcal I}_3(I,\u,\c; x)
=\int\limits_{\substack{\boldsymbol{y}\in[1,+\infty[^{10}\\
		y_1y_2y_7y_8y_{10}\leq x\\y_3y_4y_7y_9y_{10}\leq x\\y_5y_6y_8y_9y_{10}\leq x}}y_7y_8y_9y_{10}^2\,d\boldsymbol{y}
=
\int\limits_{\substack{y_8,y_9,y_{10}\in[1,+\infty[\\y_8y_9y_{10}\leq x}}y_8y_9y_{10}^2
\int\limits_{\substack{y_1,\ldots,y_7\in[1,+\infty[\\
		y_1y_2y_7\leq x/y_8y_{10}\\y_3y_4y_7\leq x/y_9y_{10}\\y_5y_6\leq x/y_8y_9y_{10}}}y_7\,d\boldsymbol{y}.$$
By using the software WX Maxima we obtain that 
\begin{eqnarray*}
{\mathcal I}_3(I,\u,\c; x)
\!\!\!&=&\!\!\!x^3\left(\frac{47}{16128}\ln^7(x)-\frac{217}{11520}\ln^6(x)+
\frac{11}{240}\ln^5(x)-
\frac{1}{32}\ln^4(x)+
\frac{4}{3}\ln^3(x)\right.\\
&&\left.-
\frac{1}{4}\ln(x)-
\frac{973}{36}\right) +O\left(x^{2+\varepsilon}\right).
\end{eqnarray*}
which implies that
$\ds K_3(c_3,\| \|_\infty):=\lim_{x\rightarrow \infty} {\mathcal I}_3(I,\u,\c; x) ~x^{-3} (\ln x)^{-7}=\frac{47}{16128}.$
\subsection{Computation of $K_n(s_n,\| \|_\infty) $ $(n=2,3)$ in corollary \ref{corollary2}}
$\bullet$ In dimension $n=2$: We have $I=\{\e_1,\e_2,\e_1+\e_2\}$, $\u=(1,1,1)$ and $\c=\zerob$. It follows that
$$\mathcal{I}_2(I,\u,\c; x)=\int\limits_{\substack{y_1,y_2,y_3\in[1,+\infty[\\y_1y_3\leq x\\y_2y_3\leq x}}\frac{d\boldsymbol{y}}{y_1y_2y_3}=\frac{1}{3}\ln^3(x).$$
which implies that
$K_2(s_2,\| \|_\infty):=\lim_{x\rightarrow \infty} {\mathcal I}_2(I,\u,\c; x) ~ (\ln x)^{-3}=\frac{1}{3}.$

$\bullet$ In dimension $n=3$: We have $I=\{\e_1,\e_2,\e_3,\e_1+\e_2,\e_1+\e_3,\e_2+\e_3,\e_1+\e_2+\e_3\}$,
$\u=\unb$ and $\c=\zerob$. It follows that
$$
\mathcal{I}_3(I,\u,\c; x)
=\int\limits_{\substack{y_1,\cdots,y_7\in[1,+\infty[\\y_1y_4y_5y_7\leq x\\y_2y_4y_6y_7\leq x\\y_3y_5y_6y_7\leq x}}\frac{d\boldsymbol{y}}{y_1y_2y_3y_4y_5y_6y_7}
=\ln^7(x)
\int\limits_{\substack{z_1,\ldots,z_7\in[0,+\infty[\\z_1+z_4+z_5+z_7\leq 1\\z_2+z_4+z_6+z_7\leq 1\\ z_3+z_5+z_6+z_7\leq 1}}d\boldsymbol{z}\\
=\frac{11}{3360}\ln^7(x). $$
We deduce that
$K_3(s_3,\| \|_\infty):=\lim_{x\rightarrow \infty} {\mathcal I}_3(I,\u,\c; x) ~ (\ln x)^{-7}=\frac{11}{3360}.$

\subsection{Computation of $K_n(u_n,\| \|_\infty) $ $(n=2,3)$ in corollary \ref{corollary3}}
$\bullet$ In dimension $n=2$: We have $I=\{\e_1,\e_2\}$, $\u=(1,1)$ and $\c=\zerob$. It follows that
$$\mathcal{I}_2(I,\u,\c; x)=\int\limits_{\substack{y_1,y_2\in[1,+\infty[\\y_1\leq x,~y_2\leq x}}\frac{d\boldsymbol{y}}{y_1y_2}=\ln^2(x)$$
and therefore that
$K_2(u_2,\| \|_\infty):=\lim_{x\rightarrow \infty} {\mathcal I}_2(I,\u,\c; x) ~ (\ln x)^{-2}=1.$

$\bullet$ In dimension $n=3$: We have $I=\{\e_1,\e_2,\e_3,\e_1+\e_2,\e_1+\e_3,\e_2+\e_3\}$, $\u=\unb$ and $\c=\zerob$. It follows that
$$
\mathcal{I}_3(I,\u,\c; x)
=\int\limits_{\substack{y_1,\cdots,y_6\in[1,+\infty[\\y_1y_4y_5\leq x\\y_2y_4y_6\leq x\\y_3y_5y_6\leq x}}\frac{d\boldsymbol{y}}{y_1y_2y_3y_4y_5y_6}
=\ln^6(x)
\int\limits_{\substack{z_1,\ldots,z_6\in[0,+\infty[\\z_1+z_4+z_5\leq 1\\z_2+z_4+z_6\leq 1\\ z_3+z_5+z_6\leq 1}}d\boldsymbol{z}
=\frac{11}{480}\ln^6(x)$$
and therefore that
$K_3(u_3,\| \|_\infty):=\lim_{x\rightarrow \infty} {\mathcal I}_3(I,\u,\c; x) ~ (\ln x)^{-6}=\frac{11}{480}.$
\subsection{Computation of $K_n(v_n,\| \|_\infty) $ $(n=2,3)$ in corollary \ref{corollary4}}
$\bullet$ In dimension $n=2$: We have $I=\{\e_1,\e_2,\e_1+\e_2\}$, $\u=(1,1,1)$ and $\c=(1,1)$. It follows that 
$\ds \mathcal{I}_2(I,\u,\c; x)=\int\limits_{\substack{y_1,y_2,y_3\in[1,+\infty[\\y_1y_3\leq x,~y_2y_3\leq x}}y_3d\boldsymbol{y}=x^2\ln(x)-\frac{3}{2}x^2+2x-\frac{1}{2}$
and therefore that 
$K_2(v_2,\| \|_\infty):=\lim_{x\rightarrow \infty} {\mathcal I}_2(I,\u,\c; x) ~x^{-2}~ (\ln x)^{-1}=1.$

$\bullet$ In dimension $n=3$: We have $I=\{\e_1,\e_2,\e_3,\e_1+\e_2,\e_1+\e_3,\e_2+\e_3,\e_1+\e_2+\e_3\}$, $\u=\unb$ and $\c=(1,1,1)$. It follows that
\begin{align*}
\mathcal{I}_3(I,\u,\c; x)
&=\int\limits_{\substack{y_1,\cdots,y_7\in[1,+\infty[\\y_1y_4y_5y_7\leq x\\y_2y_4y_6y_7\leq x\\y_3y_5y_6y_7\leq x}}y_4y_5y_6y_7^2d\boldsymbol{y}
=\ln^7(x)
\int\limits_{\substack{z_1,\ldots,z_7\in[0,+\infty[\\z_1+z_4+z_5+z_7\leq 1\\z_2+z_4+z_6+z_7\leq 1\\ z_3+z_5+z_6+z_7\leq 1}}x^{z_1+\cdots+z_3+2z_4+\cdots+2z_6+3z_7}d\boldsymbol{z}\\
&=\ln^7(x)
\int\limits_{\substack{z_5,z_6,z_7\in[0,+\infty[\\z_5+z_6+z_7\leq 1 }}\int\limits_{\substack{z_1,z_2,z_4\in[0,+\infty[\\z_1+z_4\leq 1-z_5-z_7\\z_2+z_4\leq 1-z_6-z_7}}
x^{z_1+z_2+2z_4+\cdots+2z_6+3z_7}\frac{x^{1-z_5-z_6-z_7}-1}{\ln(x)}\,dz_{1,2,4}\,dz_{5,6,7}\\
&=\ln^4(x)
\int\limits_{\substack{z_5,z_6,z_7\in[0,+\infty[\\z_5+z_6+z_7\leq 1}}
\int\limits_{z_4=0}^{1-z_7-z_5\vee z_6}
x^{2z_4+\cdots+2z_6+3z_7}(x^{1-z_5-z_6-z_7}-1)(x^{1-z_4-z_5-z_7}-1)\times\\
&\quad\times(x^{1-z_4-z_6-z_7}-1)\,dz_4\,dz_{5,6,7}
\end{align*}
By symmetry in $z_5$ and $z_6$  we get
\begin{align*}
\mathcal{I}_3(I,\u,\c; x)
&=2\ln^4(x)
\left(\int\limits_{z_7=0}^{1}
\int\limits_{z_6=0}^{(1-z_7)/2}\int\limits_{z_5=0}^{z_6}
\int\limits_{z_4=0}^{1-z_6-z_7}+
\int\limits_{z_7=0}^{1}
\int\limits_{z_6=(1-z_7)/2}^{1-z_7}\int\limits_{z_5=0}^{1-z_6-z_7}
\int\limits_{z_4=0}^{1-z_6-z_7}\right)\\
&=\frac{1}{16}x^3\ln^4(x)-
\frac{1}{4}x^3\ln^3(x)+
\frac{5}{2}x^3\ln(x)-\frac{67}{12}x^3 + O(x^{2+\varepsilon}).
\end{align*}
We deduce that
$K_3(v_3,\| \|_\infty):=\lim_{x\rightarrow \infty} {\mathcal I}_3(I, \u,\c; x) ~ x^{-3}~(\ln x)^{-4}=\frac{1}{16}.$
\subsection{Computation of $K_n(.,\| \|_d) $ $(n=2,3)$ in corollaries \ref{corollary5} and \ref{corollary6}}
\subsubsection{Sargos's volume constant}
First we will recall some notations from \S 2.3.1 of \cite{Toricmanin}.
Let $Q(\X)= \sum_{\alphab \in supp(Q)} a_{\alphab} \X^{\alphab}$ be
a generalized polynomial with  positive coefficients that depends upon all the variables $X_1,\dots,X_n$.
We apply the discussion in \cite{sargosthese} (see also \cite{sargosorsay}) to define a ``volume constant'' for $Q$.\par
By definition, the Newton polyhedron of $Q$ (at infinity)  is the set
${\mathcal E}^{\infty}(Q):= \left  (conv(supp(Q))-\R_+^n\right).$\par
Let $G_0$ be the smallest face of ${\mathcal E}^{\infty}(Q)$ which meets
the diagonal $\Delta =\R_+ \unb $. We denote by $\sigma_0$ the unique positive real number
$t$   that satisfies  $t^{-1}\unb \in G_0$. Thus, there exists a unique vector
subspace $\overrightarrow G_0$ of largest codimension $\rho_0$ such that
$G_0 \subset \sigma_0^{-1} \unb + \overrightarrow G_0.$
 Both $\rho_0, \sigma_0$ evidently depend upon $Q,$ but it is not necessary to indicate
 this in the notation. We also set $Q_{G_0}(X)=\sum_{\alphab \in G_0} a_{\alphab} \X^{\alphab}$.\par
There exist finitely many facets of ${\mathcal E}^{\infty}(Q)$ that intersect in $G_0.$ We denote their normalized polar vectors   by $\lambdab_1,\dots,\lambdab_N$. \par
By a permutation of the coordinates $X_i$ one can suppose that $\oplus_{i=1}^{\rho_0}\R \e_i
\oplus \overrightarrow{G_0}=\R^n\,,$ and that $\{\e_{m+1},\dots,\e_n\}$ is the set of vectors to which $G_0$ is parallel (i.e. for which $ G_0=G_0-\R_+ \e_i$). If $G_0$ is compact then $m = n$. \\
Set   $\Lambda =Conv\{\zerob,\lambdab_1,\dots,\lambdab_N, \e_{\rho_0+1},\dots,
\e_n\}$. It  follows that $dim \Lambda = n.$
\begin{definition}\label{sargosconstant}
The {\it volume constant} associated to $Q$ is:
$$
A_0(Q):= n!~Vol(\Lambda)~\int_{[1,+\infty)^{n-m}}\left(\int_{\R_+^{m-\rho_0}}P_{G_0}^{-\sigma_0}
(\unb,\x, \y) ~d\x\right) d\y\,.
$$
\end{definition}
 In (\cite{sargosthese}, chap 3, th. 1.6) (also see \cite{sargosorsay}),  P. Sargos proved
the following important result:

{\bf Theorem (P. Sargos):}
{\it Let $Q$ be a generalized polynomial with positive coefficients. Then
$\ds s\mapsto Y(Q;s):=\int_{[1,+\infty)^n } Q(\x)^{-s} ~d\x$ converges absolutely in $\{\Re s > \sigma_0\},$ and has a
meromorphic continuation to $\C$
with largest pole at $s = \sigma_0$ of order $\rho_0$. In addition, the volume constant of $Q$ is given by
\begin{equation}\label{sargosth}
 A_0(Q) = \lim_{s\rightarrow \sigma_0} (s-\sigma_0)^{\rho_0} ~Y(Q;s) >0.
 \end{equation}}
\subsubsection{Mellin's Formula}
We will also use the following classical Mellin's formula:\\
Let $w_1,\dots, w_r \in \C$ such that $\Re (w_i)>0$ $\forall i=0,\dots, r$,
Let $\rho_1,\dots, \rho_r >0$. Then, for $s\in \C$ verifying
$\Re(s) >\rho_1+\dots +\rho_r$, we have :
\begin{equation}\label{Mellinformula}
\frac{\Gamma (s)}{\left(\sum_{k=0}^r w_k \right)^{s}}
=
\frac{1}{(2\pi i)^r}
\int_{(\rho_1)}\dots
\int_{(\rho_r)}
 \frac{\Gamma(s-z_1-\dots-z_r)~\prod_{i=1}^r \Gamma (z_i) ~d\z}
{w_0^{s-z_1-\dots-z_r}\left(\prod_{k=1}^r w_k^{z_k}\right)},
\end{equation}
where the notation $\int_{(\rho)}$ denote the integral on the vertical line $\Re(s)=\rho$.

\subsubsection{Computation of $K_n(c_n,\| \|_d) $ $(n=2,3)$ in corollary \ref{corollary5}}
$\bullet$  In dimension $n=2$: Corollary \ref{corollary5} implies that
\begin{equation}\label{k2c2d}
K_2(c_2,\| \|_d)=\frac{d^{4}}{24}~ A_0(\T,P),
\end{equation}
 where $A_0(\T,P)$ is the mixed volume constant (see \S  2.3.3  of \cite{Toricmanin}) associated to the polynomial
 $P= X_1^{d}+ X_2^{d}$ and  the pair $\mathcal T =\left(\tilde{I}, \u=\left(u(\betab)\right)_{\betab \in\tilde{I}} \right)$, where\\
$\tilde{I}=\{\frac{1}{2}(\boldsymbol{e_1}+\boldsymbol{e_2}),\boldsymbol{e_1},\boldsymbol{e_2}\}$, and
$ u\left(\frac{1}{2}(\e_1+\e_2)\right)=1 {\mbox { and }} u(\e_1)=u(\e_2)=2.$\\
It follows then from the construction given in \S 2.3.3 of \cite{Toricmanin} that
\begin{equation}\label{k2c2dbis}
A_0(\T,P) = A_0(Q),
\end{equation}
where $A_0(Q)$ is the volume constant associated to the polynomial
$$Q(X_1,X_2,X_3,X_4,X_5)=X_1^{d/2} X_2^d X_3^d+X_1^{d/2} X_4^d X_5^d.$$
By using notations of \S 7.7.1, we have
$$\mathcal E^\infty(Q):= conv\left(supp(Q)-\R_+^5\right)=  conv\left(\{(d/2,d,d,0,0), (d/2,0,0,d,d)\}-\R_+^5\right),$$
$$G_0=\operatorname{conv}\{(d/2,d,d,0,0),(d/2,0,0,d,d)\},\quad \sigma_0 =2/d \quad {\mbox{ and }} \quad  \rho_0=4.$$

Sargos's Theorem above implies then that
\begin{equation}\label{k2c2dbis2}
A_0(Q)=\lim_{s\rightarrow 2/d} \left(s-\frac{2}{d}\right)^4 Y(Q; s).
\end{equation}
We will now compute the principal part of the integral $Y(Q;s)$.\\
First we remark that for $\Re(s)>2/d$, we have
	\begin{equation}\label{yqsc21}
	Y(Q;s)=\int_{[1,\infty)^5}(x_1^{d/2}x_2^dx_3^d+x_1^{d/2}x_4^dx_5^d)^{-s}dx_{1,2,3,4,5}
	=\frac{2}{d}\frac{1}{s-\frac{2}{d}} \displaystyle\int_{[1,\infty)^4}(x_2^dx_3^d+x_4^dx_5^d)^{-s}dx_{2,3,4,5}
 	\end{equation}
	Mellin's formula (\ref{Mellinformula}) implies that  for $\Re(s)>2/d$,
	\begin{eqnarray*}
	\displaystyle\int_{[1,\infty)^4}(x_2^dx_3^d+x_4^dx_5^d)^{-s}dx_{2,3,4,5}
	&=\dfrac{1}{2\pi i}\displaystyle\int_{[1,\infty)^4}\int_{(2/d)}\dfrac{\Gamma(s-z)\Gamma(z)}{\Gamma(s)}(x_2x_3)^{-d(s-z)}(x_4x_5)^{-dz}dx_{2,3,4,5}dz\\
	&=\dfrac{1}{2\pi i}\displaystyle\int_{(2/d)}\dfrac{\Gamma(s-z)\Gamma(z)}{d^4\Gamma(s)}\dfrac{1}{\left[(s-z)-\frac{1}{d}\right]^2}\dfrac{1}{\left[z-\frac{1}{d}\right]^2}dz
    \end{eqnarray*}
	Moving the integration line to left until $\frac{1}{2d}$ and using residues theorem imply that
	\begin{eqnarray}\label{trucbien}
	\displaystyle\int_{[1,\infty)^4}(x_2^dx_3^d+x_4^dx_5^d)^{-s}dx_{2,3,4,5}
	&=\dfrac{1}{2\pi i}\displaystyle\int_{\left(\frac{1}{2d}\right)}\dfrac{\Gamma(s-z)\Gamma(z)}{d^4\Gamma(s)}\dfrac{1}{\left[(s-z)-\frac{1}{d}\right]^2}\dfrac{1}{\left[z-\frac{1}{d}\right]^2}dz\quad\quad\quad  \\
	&+\dfrac{\Gamma^{(1)}\left(\frac{1}{d}\right)\Gamma\left(s-\frac{1}{d}\right)}{d^4\Gamma(s)\left[s-\frac{2}{d}\right]^2}
	-\dfrac{\Gamma\left(\frac{1}{d}\right)\Gamma^{(1)}\left(s-\frac{1}{d}\right)}{d^4\Gamma(s)\left[s-\frac{2}{d}\right]^2}
	+\dfrac{2\Gamma\left(\frac{1}{d}\right)\Gamma\left(s-\frac{1}{d}\right)}{d^4\Gamma(s)\left[s-\frac{2}{d}\right]^3}\nonumber
	\end{eqnarray}
	Since the integral in the  right side of (\ref{trucbien}) defines a holomorphic  function in $\Re (s)>\frac{3}{2d}$, we deduce
	by using in addition (\ref{yqsc21}) that
	\begin{equation}\label{k2c2dbis2bis}
A_0(Q)=\lim_{s\rightarrow 2/d} \left(s-\frac{2}{d}\right)^4 Y(Q; s)=
 	\dfrac{4\Gamma\left(\frac{1}{d}\right)^2}{d^5\Gamma\left(\frac{2}{d}\right)}.
 	\end{equation}
 	Combining (\ref{k2c2d}), (\ref{k2c2dbis}) and (\ref{k2c2dbis2bis}) implies that
 	\begin{equation}\label{k2c2df}
K_2(c_2,\| \|_d)=\dfrac{1}{6d}\dfrac{\Gamma\left(\frac{1}{d}\right)^2}{\Gamma\left(\frac{2}{d}\right)}.
\end{equation}

$\bullet$ In dimension $n=3$: Corollary \ref{corollary5} implies that
\begin{equation}\label{k3c3d}
K_3(c_3,\| \|_d)=\frac{d^{8}}{72 \times 7!}~ A_0(\T,P),
\end{equation}
 where $A_0(\T,P)$ is the mixed volume constant (see \S  2.3.3  of \cite{Toricmanin}) associated to the polynomial
 $P= X_1^{d}+ X_2^{d}+X_3^d$ and  the pair $\mathcal T =\left(\tilde{I}, \u=\left(u(\betab)\right)_{\betab \in\tilde{I}} \right)$, where
$$\tilde{I}=\left\{\frac{1}{3}(\boldsymbol{e_1}+\boldsymbol{e_2}+\boldsymbol{e_3}),\frac{1}{2}(\boldsymbol{e_1}+\boldsymbol{e_2}),\frac{1}{2}(\boldsymbol{e_1}+\boldsymbol{e_3}),\frac{1}{2}(\boldsymbol{e_2}+\boldsymbol{e_3}),\boldsymbol{e_1},\boldsymbol{e_2},\boldsymbol{e_3}\right\}$$
 and ${\bf u}=(1,1,1,1,2,2,2)$.
It follows then from the construction given in \S 2.3.3 of \cite{Toricmanin} that
\begin{equation}\label{k3c3dbis}
A_0(\T,P) = A_0(Q),
\end{equation}
where $A_0(Q)$ is the volume constant associated to the polynomial
$$Q(x_1,\ldots,x_{10})=x_1^{d/3}x_2^{d/2}x_3^{d/2}x_5^{d}x_8^{d}+x_1^{d/3}x_2^{d/2}x_4^{d/2}x_6^dx_9^{d}+x_1^{d/3}x_3^{d/2}x_4^{d/2}x_7^{d}x_{10}^{d}.$$
By using notations of \S 7.7.1, we have\\
$G_0=\operatorname{conv}\{(\frac{d}{3},\frac{d}{2},\frac{d}{2},0,d,d,0,0,0,0),(\frac{d}{3},\frac{d}{2},0,\frac{d}{2},0,0,d,d,0,0), (\frac{d}{3},0,\frac{d}{2},\frac{d}{2},0,0,0,0,d,d)\},
$
$$\sigma_0 = 3/d\quad {\mbox{ and }} \quad  \rho_0=8.$$

Sargos's Theorem implies then that
\begin{equation}\label{k3c3dbis2}
A_0(Q)=\lim_{s\rightarrow 3/d} \left(s-\frac{3}{d}\right)^8 Y(Q; s).
\end{equation}
We will now compute the principal part of the integral $Y(Q;s)$ at $s=3/d$.\\
Mellin's formula (\ref{Mellinformula}) implies that  for $\Re(s)>3/d$,
$$Y(Q; s)=\dfrac{24}{d^{10}\left[s-\frac{3}{d}\right]}\dfrac{1}{(2\pi i)^2}
 	\displaystyle\underset{\substack{(\rho_1)=\frac{3}{2d}\\(\rho_2)=\frac{3}{2d}}}{\int}
 	\dfrac{\Gamma(s-z_1-z_2)\Gamma(z_1)\Gamma(z_2)\Gamma(s)^{-1} \, dz_{1,2}}{
 	\left(s-z_1-z_2-\frac{1}{d}\right)^2\left(z_1+z_2-\frac{2}{d}\right)
 	\prod\limits_{j=1}^2\left[\left(z_j-\frac{1}{d}\right)^2\left(s-z_j-\frac{2}{d}\right)\right]  }.$$
By using the residue theorem, we obtain (the details of computation are left to the reader) that for $\Re (s)>3/d$,
\begin{equation}\label{19dec}
Y(Q;s) = \dfrac{372\Gamma\left(\frac{1}{d}\right)^2\Gamma\left(s-\frac{2}{d}\right)}{d^{10}\Gamma(s)\left[s-\frac{3}{d}\right]^8} +\dfrac{H(s)}{\left[s-\frac{3}{d}\right]^7},
\end{equation}
where $H$ is a holomorphic function in the bigger domain $\Omega=\left\{\Re (s) >\frac{5}{2d}\right\}$.\\
 Combining (\ref{k3c3d}), (\ref{k3c3dbis}), (\ref{k3c3dbis2}) and  (\ref{19dec}) implies that
 \begin{equation}\label{k3c3df}
K_3(c_3,\| \|_d)=\dfrac{31~\Gamma\left(\frac{1}{d}\right)^3}{30240~ d^{2}~\Gamma\left(\frac{3}{d}\right)}.
\end{equation}	
	
	\subsubsection{Computation of $K_n(v_n,\| \|_d) $ $(n=2,3)$ in corollary \ref{corollary6}}
	
$\bullet$  In dimension $n=2$: Corollary \ref{corollary6} implies that
\begin{equation}\label{k2v2d}
K_2(v_2,\| \|_d)=\frac{d^{2}}{4}~ A_0(\T,P),
\end{equation}
 where $A_0(\T,P)$ is the mixed volume constant (see \S  2.3.3  of \cite{Toricmanin}) associated tothe polynomial
 $P= X_1^{d}+ X_2^{d}$ and  the pair $\mathcal T =\left(\tilde{I}, \u=\left(u(\betab)\right)_{\betab \in\tilde{I}} \right)$, where\\
$\tilde{I}=\{\frac{1}{2}(\boldsymbol{e_1}+\boldsymbol{e_2}),\boldsymbol{e_1},\boldsymbol{e_2}\}$, and
$u\left(\frac{1}{2}(\e_1+\e_2)\right)= u(\e_1)=u(\e_2)=1.$\\
It follows then from the construction given in \S 2.3.3 of \cite{Toricmanin} that
\begin{equation}\label{k2c2dbis}
A_0(\T,P) = A_0(Q),
\end{equation}
where $A_0(Q)$ is the volume constant associated to the polynomial
$$Q(X_1,X_2,X_3)=X_1^{d/2} X_2^d +X_1^{d/2} X_3^d.$$
By using notations of \S 7.7.1, we have
$$G_0=\operatorname{conv}\{(d/2,d,0),(d/2,0,d)\},~\sigma_0 =2/d,~  \rho_0=2 {\mbox{ and }} 
\Lambda=\operatorname{conv}\left\{{\bf 0},\frac{1}{d}(2,0,0),\frac{1}{d}(0,1,1) ,\boldsymbol{e_3}\right\}.$$
It follows then from Definition \ref{sargosconstant} above that
$$A_0(\mathcal{T},P)=A_0(Q)=3!\operatorname{Vol}(\Lambda)\int_{\R_+}Q({\bf 1},x_3)^{-2/d}\,dx_3=\dfrac{2}{d^3}\dfrac{\Gamma(1/d)^2}{\Gamma(2/d)}.$$
By using in addition (\ref{k2v2d}) we obtain that
\begin{equation}\label{k2v2df}
K_2(v_2,\| \|_d) =\dfrac{1}{2d}\dfrac{\Gamma(1/d)^2}{\Gamma(2/d)}.
\end{equation}

$\bullet$ In dimension $n=3$: Corollary \ref{corollary6} implies that
\begin{equation}\label{k3v3d}
K_3(v_3,\| \|_d)=\frac{d^{5}}{3\times 4!}~ A_0(\T,P),
\end{equation}
 where $A_0(\T,P)$ is the mixed volume constant (see \S  2.3.3  of \cite{Toricmanin}) associated to the polynomial
 $P= X_1^{d}+ X_2^{d}+X_3^d$ and  the pair $\mathcal T =\left(\tilde{I}, \u=\left(u(\betab)\right)_{\betab \in\tilde{I}} \right)$, where
$$\tilde{I}=\left\{\frac{1}{3}(\boldsymbol{e_1}+\boldsymbol{e_2}+\boldsymbol{e_3}),\frac{1}{2}(\boldsymbol{e_1}+\boldsymbol{e_2}),\frac{1}{2}(\boldsymbol{e_1}+\boldsymbol{e_3}),\frac{1}{2}(\boldsymbol{e_2}+\boldsymbol{e_3}),\boldsymbol{e_1},\boldsymbol{e_2},\boldsymbol{e_3}\right\}$$
 and ${\bf u}=(1,1,1,1,1,1,1)$.
It follows then from the construction given in \S 2.3.3 of \cite{Toricmanin} that
\begin{equation}\label{k3v3dbis}
A_0(\T,P) = A_0(Q),
\end{equation}
where $A_0(Q)$ is the volume constant associated to the polynomial
$$Q(x_1,\ldots,x_{7})=x_1^{d/3}x_2^{d/2}x_3^{d/2}x_5^{d}+x_1^{d/3}x_2^{d/2}x_4^{d/2}x_6^d+x_1^{d/3}x_3^{d/2}x_4^{d/2}x_7^{d}.$$
By using notations of \S 7.7.1, we have
$$G_0=\operatorname{conv}\{(d/3,d/2,d/2,0,d,0,0),(d/3,d/2,0,d/2,0,d,0),\\(d/3,0, d/2,d/2,0,0,d)\},$$
$$\sigma_0 = 3/d\quad {\mbox{ and }} \quad  \rho_0=5.$$

Sargos's Theorem above implies then that
\begin{equation}\label{k3v3dbis2}
A_0(Q)=\lim_{s\rightarrow 3/d} \left(s-\frac{3}{d}\right)^5 Y(Q; s).
\end{equation}
Using Mellin's formula (\ref{Mellinformula})  as in the proof of (\ref{19dec}) implies that for $\Re (s)>3/d$, we have
	\begin{eqnarray*}
	Y(Q;s)=\displaystyle\int_{[1,\infty)^7}(x_1^{d/3}x_2^{d/2}x_3^{d/2}x_5^d+x_1^{d/3}x_2^{d/2}x_4^{d/2}x_6^d+x_1^{d/3}x_3^{d/2}x_4^{d/2}x_7^d)^{-s}dx_{1,2,3,4,5,6,7}\\
	= \dfrac{36\Gamma\left(s-\frac{2}{d}\right)\Gamma\left(\frac{1}{d}\right)^2}{d^7\Gamma(s)\left[s-\frac{3}{d}\right]^5}+\dfrac{H(s)}{\left[s-\frac{3}{d}\right]^4},
	\end{eqnarray*}
	where $H$ is a holomorphic function in the domain $\Omega=\left\{\Re (s) >\frac{5}{2d}\right\}$.\\
	We deduce that
	$$A_0(\mathcal{T},P)=A_0(Q)=\dfrac{36~\Gamma\left(1/d\right)^3}{d^7~\Gamma\left(3/d\right)}.$$

It follows then from (\ref{k3v3d}) that
\begin{equation}\label{k3v3df}
K_3(v_3,\| \|_d)=\dfrac{\Gamma\left(1/d\right)^3}{2 ~d^2 ~\Gamma\left(3/d\right)}.
\end{equation}

\section{Acknowledgments}
Part of this work was done when the third author visited UJM-Saint-Etienne in June 2019. He thanks for the invitation and hospitality. The research
of the third author was also financed by NKFIH in Hungary, within the framework of the 2020-4.1.1-TKP2020 3rd thematic programme of the University of P\'ecs.

\end{document}